\documentclass[12pt]{article}
\usepackage{latexsym,amsfonts,amssymb,mathrsfs,amsmath,amsthm,caption,subfigure}
\usepackage[pagewise]{lineno}
\usepackage{ulem}
\usepackage{color}
\usepackage{xcolor}
\usepackage{colordvi,multicol}
\usepackage{graphicx}
\usepackage{epstopdf}
\usepackage{bm}
\usepackage{multirow}
\usepackage{float}
\usepackage{bbding}
\renewcommand{\theequation}{\arabic{section}.\arabic{equation}}
\usepackage{indentfirst}
\usepackage{amsmath}
\usepackage{url}
\usepackage[colorlinks=true,
citecolor=blue,
urlcolor=blue,
linkcolor=blue]{hyperref}  

\numberwithin{equation}{section}
\numberwithin{figure}{section}

\newtheorem{thm}{Theorem}[section]

\newtheorem{pro}[thm]{Proposition}
\newtheorem{defi}{Definition}[section]
\newtheorem{rem}{Remark}[section]

\newtheorem{ass}{Assumption}[section]

\makeatletter
\def\@makefnmark{}
\makeatother

\begin{document}
	\title{\bf {Robust Optimal Reinsurance, Investment, and Surplus Allocation for Epstein--Zin Preferences}\\
	}
	
	\renewcommand{\thefootnote}{}
	\footnotetext{
		$^{*}$ Corresponding author.\\
		E-mail address: jyguo@nankai.edu.cn,\ jianxuanli@mail.nankai.edu.cn,\
		qianqzhou@yeah.net
	}
	\author{\small{Junyi Guo$^{a}$, Jianxuan Li$^{b}$, Qianqian Zhou$^{c,*}$}   \\
		\small {
			{
				$^a$School of Mathematical Sciences and LPMC,
				Nankai University, Tianjin, 300071, PR China}}\\
		\small {
			{
				$^b$School of Mathematical Sciences,
				Nankai University, Tianjin, 300071, PR China}}\\
		\small {
			{
				$^c$School of Mathematical Sciences,
				Tianjin University of Technology, Tianjin, 300384, PR China}}}
	
	\date{}
	\maketitle

	\vskip 0.5cm \noindent{\bf Abstract}\quad

	In this paper, we investigate the robust optimal reinsurance, investment, and internal surplus distribution (i.e., consumption) problem for an insurer with Epstein--Zin recursive preferences in an incomplete market. It is assumed that the insurer can allocate wealth to a financial market consisting of a risk-free asset and a risky asset, where the price process of the risky asset follows a diffusion process with a stochastic drift rate governed by an Ornstein-Uhlenbeck (O-U) process. 
	For both the unit and non-unit elasticity of intertemporal substitution (EIS) cases, by applying the classical dynamic programming approach, we derive explicit solutions for the optimal robust reinsurance, investment, and consumption strategies and also  verify that the obtained solutions indeed solve the optimal control problem. 
	Furthermore, we compare the robust solutions with their non-robust counterparts, and the comparative results shown in the figures are consistent with economic intuition.
	Finally, we contrast the exact solutions with the Campbell--Shiller approximation and assess the accuracy of the approximation method.

	\smallskip
	
	
	\smallskip
	
	
	\smallskip
	
	\noindent {\bf JEL Classiﬁcation }\quad G22 $\cdot$ G11

	\section{Introduction}

	In the field of actuarial science, the issue of the optimal reinsurance--investment strategy for insurance companies has always been a core topic of concern for both the academic community and the practical sector. On the one hand, insurance companies transfer part of their underwriting risks through reinsurance arrangements; on the other hand, they participate in financial market investments through asset allocation to achieve specific business goals, such as maximizing the expected utility of terminal wealth, minimizing the probability of ruin, and mean-variance optimization. A rich body of research has emerged on this topic, see, for example, Bai and Guo \cite{BG2008}, Liang et al. \cite{LYG2011}, Zeng and Li \cite{ZL2011}, Yi et al. \cite{YLVZ2013}, Zhang et al. \cite{ZMZ2016}, Li et al. \cite{LiZengYang2018}, Chen and Yang \cite{CY2020}, and Yuan et al. \cite{YHLY2023}.
	However, the classical studies on insurer decision-making primarily rely on time-additive expected utility theory, which simplifies preference modeling but suffers from a critical limitation: it conflates two distinct behavioral traits -- risk aversion (attitude toward uncertain outcomes across states) and elasticity of intertemporal substitution (EIS, willingness to smooth consumption over time). This entanglement contradicts empirical evidence: for instance, insurers may exhibit high risk aversion toward catastrophic underwriting risks while maintaining a moderate EIS to optimize long-term investment returns. Therefore, models based on time-additive expected utility theory fail to capture the nuanced risk-return trade-offs in dynamic insurer decision-making.

	The stochastic differential utility (SDU), also known as the Epstein--Zin recursive utility and proposed by Duffie and Epstein \cite{DE1992} as a natural continuous-time extension of the discrete-time recursive preferences developed by Epstein and Zin \cite{EZ1989}, addresses this shortcoming by disentangling risk aversion and EIS, enabling a more realistic characterization of dynamic preferences.
	This utility specification has been widely adopted in optimal consumption--portfolio problems, asset pricing, retirement planning, and longevity risk management. 
	For instance, by applying the dynamic programming principle, Kraft et al. \cite{KSS2013} and Kraft et al. \cite{KSS2017} studied the optimal consumption--portfolio decision of an investor with recursive preferences of Epstein--Zin type.
	Matoussi and Xing \cite{MX2018} investigated a continuous-time consumption and investment problem with incomplete markets and Epstein--Zin stochastic differential utilities by using a dual problem.
	Pu and Zhang \cite{PZ2021} extended the verification theorem proposed by Kraft et al. \cite{KSS2013} to a game framework involving two control variables and studied robust investment--consumption decisions.
	Xing \cite{Xing2017} studied an optimal consumption--investment problem with an unbounded market price of risk for an agent with Epstein--Zin recursive utility, using the comparison principle for backward stochastic differential equations (BSDEs). 
	Feng et al. \cite{FTZ2026} investigated the consumption--investment problem in an incomplete market and derived optimal strategies by using the martingale optimality principle and quadratic BSDEs.
	Aase \cite{Aase2016} examined optimal consumption and pension insurance decisions over an individual's lifecycle.
	Han and Hung \cite{HH2017} investigated optimal portfolio choice, consumption, and life insurance strategies for wage earners exposed to interest rate and inflation risks.
	Li et al. \cite{LSS2025} extended Yaari's lifecycle model by incorporating mortality shock analysis into the recursive utility framework, aiming to explore the impact of mortality model uncertainty on the observed low demand for annuities. 
	Peng and Xv \cite{PX2025} focused on how wage earners can make robust optimal consumption, investment, and life insurance purchase decisions when both the price process of risky assets and the stochastic factor process exhibit uncertainty.

	However, the aforementioned studies primarily concentrate on the consumption--investment decisions of individual consumers. In contrast, within the field of property reinsurance, where the insurance company itself serves as the decision-making entity, there is a notable scarcity of literature exploring the application of the Epstein--Zin utility framework.
	As is well known, the business goals of insurance companies are not limited to external risk management (reinsurance) and investment appreciation; their ``intrinsic consumption'' of internal economic surplus -- the corporate counterpart to individual consumption -- is equally crucial. Such consumption strategies are mainly reflected in capital returns to shareholders (such as dividend distributions, share repurchases), as well as performance-based compensation distribution to employees. Reasonable ``intrinsic consumption'' not only ensures long-term stable returns for shareholders but also effectively recognizes the past contributions of employees, helping to enhance internal incentives and promoting the long-term, stable, and healthy development of insurance companies. Furthermore, from a decision-making perspective, the ``consumption'' of an insurance company and the ``consumption'' of an individual are mathematically analogous in nature. Therefore, incorporating consumption strategies into the research framework of insurance company decision-making has significant practical significance and theoretical importance.

	This study aims to develop a robust optimal reinsurance--investment--consumption framework for an insurer with Epstein--Zin preferences and then investigate:
	\begin{enumerate}
		\item How do disentangled risk aversion and EIS influence an insurer's optimal allocation of surplus to reinsurance, investment, and consumption?
		\item How does model uncertainty (ambiguity aversion) affect these strategies?
		\item Can closed-form solutions for the optimal strategies be derived under realistic assumptions about surplus dynamics and market structures?​
	\end{enumerate}

	To the best of our knowledge, the only preliminary exploration of this problem is found in Dadzie et al. \cite{DKN2025arXiv}, who addressed it by solving the associated coupled forward-backward stochastic differential equations (FBSDEs).
	In contrast, this paper adopts the dynamic programming method to solve the problem, thereby extending the theoretical framework developed by Pu and Zhang \cite{PZ2021} to the context of optimal reinsurance, investment, and consumption decisions for insurance companies.
	Given the complexity of the financial and insurance market environment, an insurance company's decision-making process is often subject to multiple sources of uncertainty. Therefore, incorporating model uncertainty into the modeling framework is of significant practical importance. This paper adopts the robust control framework proposed by Hansen and Sargent \cite{HS2001} to introduce model uncertainty into the insurer's preference specification. The aforementioned studies, such as Pu and Zhang \cite{PZ2021} and Peng and Xv \cite{PX2025}, have already investigated related issues within this framework.
	Furthermore, to more accurately capture the empirical characteristics of risky assets, including the bull and bear market patterns observed in equity markets, this paper models the price process of the risky asset as a diffusion process with a stochastic drift driven by an Ornstein-Uhlenbeck (O-U) process. This specification has also been adopted in Liang et al. \cite{LYG2011}, Guan and Liang \cite{GL2014}, and Dong et al. \cite{DRZ2023}.

	Compared with existing literature, the main innovations of this paper are as follows:
	
	First, this paper focuses on the insurance company as the primary decision-maker, integrates Epstein--Zin utility into the insurer's decision-making model, incorporates investment, reinsurance, and consumption decisions, and investigates the robust optimal decision-making problem. In contrast to Dadzie et al. \cite{DKN2025arXiv}, who adopted a risky asset price model with a constant drift rate, this paper assumes that the risky asset price follows a diffusion process with a stochastic drift rate driven by an O-U process, thereby enhancing the model's alignment with real-world financial scenarios. In terms of solution methodology, this paper employs the dynamic programming principle to solve the Hamilton-Jacobi-Bellman-Isaacs (HJBI) equation and obtains explicit closed-form solutions for the optimal strategies and the corresponding value function, which differs from the approach of Dadzie et al. \cite{DKN2025arXiv} using FBSDEs to address the problem.
	
	Second, due to the introduction of the stochastic drift rate driven by an O-U process, the comparison principle used in Pu and Zhang \cite{PZ2021} and other related studies is no longer applicable. To address this methodological challenge, this paper develops a new verification framework by exploring the relevant properties of the O-U process and adopts an auxiliary construction method. In terms of solution outcomes, in contrast to Han and Hung \cite{HH2017}, who provided only a single form of solution, this paper, in the case of non-unit EIS case, not only derives approximate solutions suitable for empirical applications using the Campbell--Shiller approximation method but also provides theoretically rigorous explicit closed-form solutions, thereby enhancing the applicability of the model.
	
	Third, the robust decision-making framework established in this paper exhibits strong generality and compatibility. As the intensity of model uncertainty approaches zero, the results naturally degenerate into the classical non-robust results derived in Kraft et al. \cite{KSS2013}.

	The rest of the paper proceeds as follows.
	Section \ref{sec2} describes the verification theorem employed in this study along with the relevant underlying assumptions. Section \ref{sec3} constructs the dynamic wealth model for the insurance company, within the framework of Epstein--Zin recursive utility, and formulates the optimal reinsurance, investment, and consumption (i.e., internal surplus distribution) problem incorporating model uncertainty. Section \ref{sec4} provides explicit closed-form solutions for the optimal strategies under both the non-unit EIS and unit EIS cases. In addition, it rigorously proves that the solutions derived for the non-unit EIS case satisfy the corresponding verification theorem. Section \ref{sec5}, by relaxing certain technical conditions required for the verification theorem in Section \ref{sec4}, further presents the Campbell--Shiller approximate solutions for the non-unit EIS case. Section \ref{sec6} concludes the paper. All proofs are relegated to the Appendix.

	\section{Assumptions and verification theorem}\label{sec2}

	In this section, we revisit the verification theorem for two-player stochastic games with a non-Lipschitz aggregator established by Pu and Zhang \cite{PZ2021}. Although the theorem has already been proven in the aforementioned work, we systematically outline its basic framework and key assumptions here to facilitate the step-by-step verification of its applicability conditions in the subsequent analysis, thereby laying a clear logical foundation for the rest of this paper.

	We begin by establishing the general framework. Let \((\Omega, \mathcal{F}, \mathbb{P})\) be a complete probability space equipped with a fixed terminal time \(T < \infty\), and let \(\{B_s\}_{s \in [0,T]}\) denote a \(d\)-dimensional Brownian motion defined on this space. Denote by \((\mathcal{F}_s)_{0 \leq s \leq T}\) the natural filtration generated by \(\{B_s\}_{s \in [0,T]}\), augmented by all \(\mathbb{P}\)-null sets of \(\mathcal{F}\) so that \(\mathcal{F}_0\) contains all such sets.
	
	Consider the state process \(X = \{X_s\}_{s \in [t,T]}\) taking values in \(\Xi \subseteq \mathbb{R}^n\), governed by the stochastic differential equation: 
	\begin{equation}\label{201}
		\begin{cases}
			\textnormal{d}X_s = \mu(s, X_s, u_s, v_s) \textnormal{d}s + \sigma(s, X_s, u_s, v_s) \textnormal{d}B_s, \\
			X_t = x,
		\end{cases}
	\end{equation}
	where \(\mu : [0, T] \times \Xi \times U \times V \to \mathbb{R}^n\) and \(\sigma : [0, T] \times \Xi \times U \times V \to \mathbb{R}^{n \times d}\) are measurable functions. Here, \(U \subseteq \mathbb{R}^m\) and \(V \subseteq \mathbb{R}^l\) denote the control spaces, \(t \in [0, T]\) is the initial time, and \(x \in \Xi\) is the initial state. The processes \(u\) and \(v\) are control strategies belonging to the admissible control sets, which will be specified in the following definition.

	\begin{defi}[\textbf{Admimissible control}]
		For \(t \in [0, T]\), let \(\mathcal{U}_t\) denote the \(t\)-admissible control set of \((\mathcal{F}_t)_{0 \leq t \leq T}\)-adapted feedback strategies \(u = \{u_s\}_{s \in [t,T]} = \{u(s, X_s)\}_{s \in [t,T]}\) taking values in \(U\), and \(\mathcal{V}_t\) denote the set of all \((\mathcal{F}_t)_{0 \leq t \leq T}\)-adapted processes \(v = \{v_s\}_{s \in [t,T]}\) taking values in \(V\).
	\end{defi}

	A control \(u \in \mathcal{U}_t\) (respectively, \(v \in \mathcal{V}_t\)) represents the strategy of the maximizing (resp., minimizing) player. The corresponding payoff functional of the game is defined as
	\begin{equation}\label{202}
		J(s, x; u, v) = J_s = \mathbb{E}_s \left[ \int_s^T f(r, X_r, u_r, v_r, J_r) \, \textnormal{d}r + \psi(X_T) \right], \quad s \in [t, T],
	\end{equation}
	where \( f : [0, T] \times \Xi \times U \times V \times \mathbb{R} \to \mathbb{R} \) is a measurable intertemporal aggregator, \( \psi : \Xi \to \mathbb{R} \) denotes the terminal cost, and \( \mathbb{E}_s \) represents the conditional expectation given \( \mathcal{F}_s \).

	The stochastic game problem of our interest is formulated as
	\begin{equation}\label{203}
		\text{(P)} \quad \sup_{u \in \mathcal{U}_t} \inf_{v \in \mathcal{V}_t} \mathbb{E} \left[ \int_t^T f(s, X_s, u_s, v_s, J_s) \, \textnormal{d}s + \psi(X_T) \right],
	\end{equation}
	and the associated value function is formally given by
	\begin{equation*}
		W(t, x) = \sup_{u \in \mathcal{U}_t} \inf_{v \in \mathcal{V}_t} J(t, x; u, v).
	\end{equation*}

	Before presenting the verification theorem, we first enumerate the conditions that must be satisfied.
	
	\begin{ass}\label{A1}
		For any \( u \in \mathcal{U}_t \) and \( v \in \mathcal{V}_t \), \eqref{201}-\eqref{202} have a unique strong solution.
	\end{ass}

	For \( w \in C^{1,2}([0, T] \times \Xi) \), define
	\begin{equation*}
		\mathcal{L}^{u,v}(w)(s,x) = w_s(s,x) + \langle \mu(s,x,u,v), w_x(s,x) \rangle + \frac{1}{2} \operatorname{tr}\left( \sigma(s,x,u,v)^\top w_{xx}(s,x) \sigma(s,x,u,v) \right).
	\end{equation*}

	We suppose that
	\begin{ass}\label{A2}
		For any \((s,x) \in [t,T] \times \Xi\) and \( u \in U, v \in V \), there exist measurable functions \( v^* : [t,T] \times \Xi \times U \to V \) and \( u^* : [t,T] \times \Xi \to U \) such that
		\begin{equation*}
			\begin{aligned}
				& v^* \in \arg \min_{v \in V} \left\{ \mathcal{L}^{u,v}(w)(s,x) + f(s,x,u,v,w(s,x)) \right\}, \\
				& u^* \in \arg \max_{u \in U} \left\{ \mathcal{L}^{u,v^*}(w)(s,x) + f(s,x,u,v^*,w(s,x)) \right\}.
			\end{aligned}
		\end{equation*}
	\end{ass}
	
	From Assumption \ref{A2}, it is clear that
	\begin{equation*}
		\mathcal{L}^{u^*,v^*}(w)(s,x) + f(s,x,u^*,v^*,w(s,x)) = \sup_{u \in U} \inf_{v \in V} \left\{ \mathcal{L}^{u,v}(w)(s,x) + f(s,x,u,v,w(s,x)) \right\}.
	\end{equation*}

	\begin{ass}\label{A3}
		There exists \( C > 0 \) such that for any \( s,x,u,v \in [t,T] \times \Xi \times U \times V \),
		
		\begin{equation*}
			f(s,x,u,v,z_1) - f(s,x,u,v,z_2) \leq C(z_1 - z_2), \quad \text{if}\ z_1 \geq z_2.
		\end{equation*}
	\end{ass}

	\begin{ass}\label{A4}
		If \( w \in C^{1,2}([0, T] \times \Xi) \cap C^{0,0}([0, T] \times \Xi) \) is a solution of the following dynamic programming equation
		\begin{equation}\label{DPE}
			\begin{cases}
				w_t(t, x) + \sup\limits_{u \in U} \inf\limits_{v \in V} \Big\{ \langle b(t, x; u, v), w_x(t, x) \rangle + \dfrac{1}{2} \operatorname{tr}\left( \sigma^\top(t, x, u, v) w_{xx}(t, x) \sigma(t, x, u, v) \right) \\
				\quad + f(t, x, u, v, w(t, x)) \Big\} = 0, \quad \text{in } [0, T] \times \Xi, \\
				w(T, x) = \psi(x), \quad x \in \Xi,
			\end{cases}
		\end{equation}
		the local martingale
		\begin{equation*}
			\int_t^T w_s(s, X_s)^\top \sigma(s, X_s, u_s, v_s) \textnormal{d}B_s
		\end{equation*}
		is a true martingale for any \( u \in \mathcal{U}_t \) and \( v \in \mathcal{V}_t \).
	\end{ass}

	\begin{ass}\label{A5}
		For any \( u \in \mathcal{U}_t \) and \( v \in \mathcal{V}_t \), if \((X_s, J_s)_{s \leq t}\) is the solution of \eqref{201}-\eqref{202}, we have \(\mathbb{E}\left[ \int_t^T |w(s, X_s) - J_s| \textnormal{d}s \right] < \infty\).
	\end{ass}

	Then we have the following verification theorem proposed by Pu and Zhang \cite{PZ2021}.
	
	\begin{thm}[\textbf{Verification theorem}]\label{VT}
		Under Assumptions~\ref{A1}--\ref{A5}, if the HJBI equation~\eqref{DPE} admits a classical solution \( w \), then the pair \( (u^*, v^*) \) constitutes an optimal strategy in \( \mathcal{U}_t \) and \( \mathcal{V}_t \), respectively. Moreover, \( w \) coincides with the value function of the max-min problem~\eqref{203}, and the following minimax identity holds:
		\begin{equation*}
			w(t, x) = J(t, x; u^*, v^*) = \sup_{u \in \mathcal{U}_t} \inf_{v \in \mathcal{V}_t} J(t, x; u, v) = \inf_{v \in \mathcal{V}_t} \sup_{u \in \mathcal{U}_t} J(t, x; u, v).
		\end{equation*}
	\end{thm}

	\section{Model formulation and robust optimization }\label{sec3}	
	
	In this section, we introduce a model for the insurer's wealth and formulate the robust optimization problem of interest.  We consider a continuous-time insurance financial market where all assets are traded continuously, and no transaction costs or taxes are imposed. Let \((\Omega, \mathcal{F}, \mathbb{P})\) be a complete probability space equipped with a filtration \(\{\mathcal{F}_t\}_{t \in [0, T]}\), where \(T > 0\) denotes a fixed finite terminal time, and $\mathbb{P}$  is a reference probability measure. Any decisions made at time \(t\) are contingent on the information set \(\mathcal{F}_t\), which represents the information available up to that time.

	\subsection{Model formulation}
	
	Suppose that the surplus process of the insurer is described by the following Cram$\acute{e}$r-Lundberg model:
	\begin{equation*}
		U_t = u + bt - S_t = u + bt - \sum_{i=1}^{N_t} Y_i,
	\end{equation*}
	where $u$ is the initial surplus,
	$b$ is the premium rate, $S_t$ representing the aggregate claims up to time $t$ is a generalized compound Poisson process, the claim-number process $N_t$ is a homogeneous Poisson process with deterministic  intensity $\lambda,$ possessing explicitly specified finite-dimensional distributions, and $\left\{ Y_i, i \geq 1 \right\}$ is a sequence of positive, independent and identically distributed (i.i.d.)  claim-amount random variables following a deterministic distribution $F(y)= \mathbb{P}(Y \le y)$.
	Here we denote the mean and second moment of $Y_i$ by $\mu_1:= \mathbb{E} [Y_i]$  and $\mu_2:= \mathbb{E} [Y^2_i].$

	Reinsurance is a core risk management tool for insurers to mitigate their exposure to extreme claim risks. At time $t$, the insurer may purchase proportional reinsurance or expand its business by underwriting new policies. Let $q(t) \in [0, +\infty)$ denote the insurer's retained risk proportion at time $t$: the insurer is liable for $100 q(t)\%$ of each claim amount, while the reinsurer assumes the remaining $100 (1-q(t))\%$.
	Correspondingly, the insurer is required to pay a portion of the premium to the reinsurer. In this paper, we assume that
	the reinsurance premium rate $\delta(q(t))$ is determined according to the expected value premium principle:
	$$
	\delta(q(t)) = (1-q(t)) (1+\theta_1)\lambda \mu_1,
	$$
	where $\theta_1$ is the safety loading of the reinsurer. Furthermore, the model incorporates the standard condition that this rate strictly dominates the pure risk premium, i.e.,
	$
	\lambda \mu_1 < b < \delta(0) = (1 + \theta_1) \lambda \mu_1,
	$
	which is indispensable for ruling out trivial or boundary solutions.

	We now turn to a brief description of the financial market, which  consists of two continuously tradable assets. One is a risk-free asset, whose price process $R(t)$ follows:
	\begin{equation}\label{301}
		\textnormal{d}R(t) = r R(t)\textnormal{d}t ,
	\end{equation}
	where $r > 0$ is the risk-free interest rate.
	The other is a risky asset, whose dynamics are inspired by the observations of Rishel \cite{Rishel1999} and the framework of Dong et al. \cite{DRZ2023}, and are modeled as follows:
	\begin{equation}\label{302}
		\textnormal{d}P(t) = P(t)\big( (a+ \sigma m(t)) \textnormal{d}t + \sigma \textnormal{d}W_{1,t} \big),
	\end{equation}
	and $m(t)$ follows an O-U process:
	\begin{equation}\label{303}
		\textnormal{d} m(t) = -\alpha m(t) \textnormal{d}t + \beta \left( \rho_1 \textnormal{d}W_{1,t} + \sqrt{1-\rho_1^2} \textnormal{d}W_{2,t} \right), \quad m(0) = m,
	\end{equation}
	where the volatility parameters $\sigma > 0$ and $\beta > 0$, the mean-reversion parameter $\alpha > 0$, and the interest rate $a > r$ are all constants.
	The correlation $|\rho_1| \leq 1$ is a constant. As long as $|\rho_1| < 1$, the market is incomplete.
	The processes $W_{1,t}$ and $W_{2,t}$ are mutually independent Brownian motions.
	The allocation of funds by insurers in the financial market is an inherent requirement for maximizing the enterprise value, and it also takes advantage of the long-term and stable nature of insurance funds. Therefore, in this paper, we assume that the insurer invests a portion of its surplus wealth in the risky asset, denoted by  $\pi(t)$ at time $t$, while investing the remaining amount in the risk-free asset.

	Apart from reinsurance planning and investment portfolio allocation, the ``intrinsic consumption'' strategy of insurance company, such as dividend distribution to shareholders and bonus allocation to employees, is also an essential component of the insurer's earnings management. We denote this consumption strategy as $c(t)$.
	
	In accordance with \eqref{301} and \eqref{302}, the evolution of the non-robust insurer's wealth $\{X_t^1\}_{t \in [0, T]}$, under a chosen strategy $(\pi, q, c)$, is governed by:
	\begin{equation}\label{304}
		\begin{aligned}
			\textnormal{d} X_t^1 & = (X_t^1 - \pi(t))\frac{\textnormal{d}R(t)}{R(t)} + \pi(t)\frac{\textnormal{d}P(t)}{P(t)} + (b-\delta(q(t))) \textnormal{d}t - q(t) \textnormal{d}\sum_{i=1}^{N_t} Y_i -c(t) \textnormal{d}t \\
			& = [\pi(t) (\sigma m(t) + a - r) + r X_t^1 + b - \delta(q(t)) -c(t)] \textnormal{d}t + \sigma \pi(t) \textnormal{d}W_{1,t} - q(t) \textnormal{d}\sum_{i=1}^{N_t} Y_i.
		\end{aligned}
	\end{equation}

	In terms of Iglehart \cite{I1969} and Grandell \cite{G1991}, we assume that the aggregate claim process $S_t$ can be approximated by the diffusion process
	\begin{equation*}\label{}
		\textnormal{d}\sum_{i=1}^{N_t} Y_i = \lambda \mu_1 \textnormal{d}t - \sqrt{\lambda \mu_2} \textnormal{d}W_{3,t},
	\end{equation*}
	where $W_{3,t}$ is a standard Brownian motion, independent of both $W_{1,t}$ and $W_{2,t}$.
	Then the dynamics of $X_t^1$ in \eqref{304} can be approximated by
	\begin{equation*}\label{305}
		\begin{aligned}
			\textnormal{d} X_t^1
			& = [\pi(t) (\sigma m(t) + a - r) + r X_t^1 + b - (1+\theta_1)\lambda \mu_1 + \lambda \mu_1 \theta_1 q(t) - c(t)] \textnormal{d}t \\
			& \quad + \sigma \pi(t) \textnormal{d}W_{1,t} + \sqrt{\lambda \mu_2} q(t) \textnormal{d}W_{3,t},
		\end{aligned}
	\end{equation*}
	where the initial wealth is $X_0^1 = x^1$.

	In line with the framework of Bodie et al. \cite{Bodie1992}, and Ma et al. \cite{MLC2023}, we define the insurer's total wealth $\{X_t\}_{t \in [0, T]}$ as
	\begin{equation*}\label{306}
		X_t = X_t^1 + (b - (1+\theta_1)\lambda \mu_1) \int_{t}^{T} e^{-r(s-t)} \textnormal{d}s,
	\end{equation*}
	where $\int_t^T e^{-r(s-t)} \textnormal{d}s$ captures the time value of continuous risk exposure over the period from $t$ to $T$. The coefficient $ b - (1+\theta_1)\lambda\mu_1 < 0$ reflects the deviation of reinsurance pricing from the fair price of full reinsurance. The rationale for this specification is as follows: by obtaining reinsurance coverage at a marginal cost $b$ that is below the full-coverage premium rate, the insurer appears to achieve a ``premium saving'', but this saving comes at the expense of retaining partial risk. According to the no-arbitrage pricing principle, such risk bearing must be reflected in wealth measurement --- future premium differentials should be converted into a current-period wealth adjustment via the discount factor.
	Then we derive the transformed wealth process as follows:
	\begin{equation}\label{307}
		\textnormal{d} X_t
		= [\pi(t) (\sigma m(t) + a - r) + r X_t + \lambda \mu_1 \theta_1 q(t) - c(t)] \textnormal{d}t
		+ \sigma \pi(t) \textnormal{d}W_{1,t} + \sqrt{\lambda \mu_2} q(t) \textnormal{d}W_{3,t},
	\end{equation}
	where the initial wealth is $X_0 := x = x^1 +  (b - (1+\theta_1)\lambda \mu_1) \int_{0}^{T} e^{-r(s-t)} \textnormal{d}s$.
	The remainder of this paper will focus on the analysis of the wealth process $X_t$  given in   \eqref{307}.

	We assume that the utility of the non-robust insurer is characterized by the continuous-time Epstein--Zin SDU, which can be represented by the following BSDE:
	\begin{equation*}
		\textnormal{d} V_s = -f(c_s, V_s) \textnormal{d} s + Z_s \textnormal{d} B_s, \quad V_T = \psi(X_T), 
	\end{equation*}
	where the Epstein--Zin aggregator is
	\begin{equation}\label{nonunitEIS}
		f(c, v) = \left(1-\frac{1}{\varphi}\right)^{-1} \delta (1-\gamma) v \left[ \left( \frac{c}{((1 - \gamma)v)^{\frac{1}{1-\gamma}}} \right)^{1-\frac{1}{\varphi}} - 1 \right]. 
	\end{equation}
	Here, \(\delta > 0\) is the rate of time preference, \(0 < \gamma \neq 1\) is the coefficient of relative risk aversion, \(0 < \varphi \neq 1\) is EIS, and \(\psi(x) :=  \frac{1}{1-\gamma}x^{1-\gamma}\).
	The case \( \varphi = 1 \) of unit EIS corresponds to the specification
	\begin{equation}\label{unitEIS}
		f(c, v) = \delta(1 - \gamma)v \left[ \ln c - \frac{1}{1 - \gamma} \ln((1 - \gamma)v) \right].
	\end{equation}
	We note that the value of $V_t$, denoted by $\mathbb{V}$, depends on the parameter $\gamma$: $\mathbb{V} = (0, \infty)$ if $\gamma < 1$, and $\mathbb{V} = (-\infty, 0)$ if $\gamma > 1$.
	
	In this study, we consider that \((\pi, q, c)\) as feedback controls contingent on the state \((x, m)\). In the absence of model uncertainty, the non-robust insurer's problem is formulated as
	\begin{equation}\label{P0}
		\text{(P0)} \quad \max_{(\pi, q, c) \in \mathcal{U}_t^0(x, m)} \mathbb{E}_t \left[ \int_t^T f(c_s, V_s) \textnormal{d}s + \psi(X_T) \right],
	\end{equation}
	where the \(t\)-admissible control set $\mathcal{U}_t^0(x, m)$, with \(x \in \mathbb{R}^+\), \(m \in \mathbb{R}\), consists of \((\mathcal{F}_t)_{0 \leq t \leq T}\)-adapted feedback strategies \(\pi(s) = \{\pi(s, X_s, m(s))\}_{s \in [t, T]}\), \(q(s) = \{q(s, X_s, m(s))\}_{s \in [t, T]}\) and \(c(s) = \{c(s, X_s, m(s))\}_{s \in [t, T]}\),  taking values in \(\mathbb{R}\), \(\mathbb{R}^+\) and  \(\mathbb{R}^+\), respectively.

	\subsection{Robust optimization problem}
	
	In this subsection, we investigate the decision-making problem of a robust insurer who possesses a relatively reliable reference model (also referred to as the baseline model). However, due to concerns about potential misspecification of the model parameters, the insurer suspects that the state dynamics described by equations \eqref{303} and \eqref{304} may be an inaccurate approximation. Consequently, the insurer is motivated to consider a set of alternative models in order to formulate a corresponding robust optimal strategy.

	To define alternative models, we consider other probability measures that are equivalent to the reference measure  $\mathbb{P}$,   defined via the Radon–Nikodym derivative as follows:
	\begin{equation*}
		\begin{split}
			\left. \frac{\textnormal{d} \mathbb{Q}^\xi}{\textnormal{d} \mathbb{P}} \right|_{\mathcal{F}_T}
			= \exp \Bigg\{ &- \int_0^T \xi_1(t) \textnormal{d}W_{1,t} - \int_0^T \xi_2(t) \textnormal{d}W_{2,t} - \int_0^T \xi_3(t) \textnormal{d}W_{3,t} \\
			&- \frac{1}{2} \int_0^T \left( \xi_1^2(t) + \xi_2^2(t) + \xi_3^2(t) \right) \textnormal{d}t \Bigg\},
		\end{split}
	\end{equation*}
	where the \(\mathbb{R}^3\)-valued process \(\xi(t) := (\xi_1(t), \xi_2(t), \xi_3(t))\), \(t \in [0, T]\), is the \textit{distortion process}. This process reflects the insurer's adjustment to the drift terms of the Brownian motions in the original model and is assumed to satisfy the necessary measurability and integrability conditions (e.g., the Novikov condition) to ensure the validity of the subsequent measure transformation.

	Then, under the measure \(\mathbb{Q}^\xi\), the processes
	\begin{equation*}
		W_{i,t}^{\mathbb{Q}^\xi} = W_{i,t} + \int_0^t \xi_i(s) \textnormal{d}s, \quad i=1,2,3,
	\end{equation*}
	is a standard Brownian motion.
	Furthermore, under the measure \( \mathbb{Q}^\xi \), the O-U process in \eqref{303} transforms into
	\begin{equation}\label{ROU}
		\begin{aligned}
			\textnormal{d} m(t)
			& = -\alpha m(t) \textnormal{d}t + \beta \left[ \rho_1 \left( \textnormal{d}W_{1,t}^{\mathbb{Q}^\xi} - \xi_1(t) \textnormal{d}t \right) + \sqrt{1-\rho_1^2} \left( \textnormal{d}W_{2,t}^{\mathbb{Q}^\xi} - \xi_2(t) \textnormal{d}t \right) \right] \\
			& = - \left[ \alpha m(t) + \beta \rho_1 \xi_1(t) + \beta \sqrt{1-\rho_1^2} \xi_2(t) \right] \textnormal{d}t + \beta \left( \rho_1 \textnormal{d}W_{1,t}^{\mathbb{Q}^\xi} + \sqrt{1-\rho_1^2} \textnormal{d}W_{2,t}^{\mathbb{Q}^\xi} \right),
		\end{aligned}
	\end{equation}
	and the wealth process \( X_t \) given in \eqref{307} evolves as
	\begin{equation}\label{RX}
		\begin{aligned}
			\textnormal{d} X_t
			& = [\pi(t) (\sigma m(t) + a - r) + r X_t + \lambda \mu_1 \theta_1 q(t) - c(t) - \pi(t) \sigma \xi_1(t) - q(t) \sqrt{\lambda \mu_2} \xi_3(t) ] \textnormal{d}t \\
			& \quad + \sigma \pi(t) \textnormal{d}W_{1,t}^{\mathbb{Q}^\xi} + \sqrt{\lambda \mu_2} q(t) \textnormal{d}W_{3,t}^{\mathbb{Q}^\xi}.
		\end{aligned}
	\end{equation}

	We assume that the insurer seeks a robust reinsurance--investment--consumption strategy that represents the best choice under the worst-case model.
	In light of this, the insurer requires a decision rule that remains effective across a range of distortion scenarios, particularly for small values of \( \xi \). To this end, we introduce a penalty term \( \frac{1}{2\Psi} \left( \xi_1^2 + \xi_2^2 + \xi_3^2 \right) \) in \eqref{P0}, where \( \Psi \) measures the strength of the preference for robustness and is a nonnegative parameter.
	Consequently, the objective function is modified as follows:
	\begin{equation}\label{PR}
		\text{(PR)} \quad \max_{(\pi, q, c) \in \mathcal{U}_t(x, m)} \min_{\xi \in \mathcal{V}_t} J(t, (x, m); (\pi, q, c), \xi),
	\end{equation}
	where for \( s \in [t, T] \),	
	\begin{equation*}\label{RJ}
		\begin{aligned}
			& J(s, (x, m); (\pi, q, c), \xi) := J_s \\
			& = \mathbb{E}_s^{\mathbb{Q}^\xi} \left[ \int_s^T \left( f(c_u, J_u) + \frac{1}{2\Psi_u} \left( \xi_1^2(u) + \xi_2^2(u) + \xi_3^2(u) \right) \right) \textnormal{d}u + \psi(X_T) \right].
		\end{aligned}
	\end{equation*}
	Here, \(\mathcal{U}_t(x, m)\) denotes the \(t\)-admissible control set, with \(x \in \mathbb{R}^+\) and \(m \in \mathbb{R}\). This set consists of \( \mathbb{F}^{W^\mathbb{Q}} \)-adapted feedback strategies:
	$$
	\pi(s) = \{\pi(s, X_s, m(s))\}_{s \in [t, T]},
	$$
	$$
	q(s) = \{q(s, X_s, m(s))\}_{s \in [t, T]},
	$$
	$$
	c(s) = \{c(s, X_s, m(s))\}_{s \in [t, T]},
	$$
	taking values in \(\mathbb{R}\), \(\mathbb{R}^+\), and \(\mathbb{R}^+\), respectively.
	\(\mathcal{V}_t\) denotes the space of all \(\mathbb{F}^{W^\mathbb{Q}}\)-adapted processes \(\xi = \{\xi_s\}_{s \in [t,T]}\) taking values in \(\mathbb{R}^3\), where \(\xi_s = (\xi_1(s), \xi_2(s), \xi_3(s))^\top\).
	The filtration \(\mathbb{F}^{W^\mathbb{Q}} = (\mathcal{F}_s^{W^\mathbb{Q}})_{0 \leq s \leq T}\) is generated by the three-dimensional Brownian motion \(W^\mathbb{Q} = \{(W_{1,s}^{\mathbb{Q}^\xi}, W_{2,s}^{\mathbb{Q}^\xi}, W_{3,s}^{\mathbb{Q}^\xi} )^\top\}_{s \in [0,T]} \).

	\begin{rem}
		The term $\frac{1}{2 \Psi} \left( \xi_1^2 + \xi_2^2 + \xi_3^2 \right)$ in the objective function is essentially a discounted relative entropy that quantifies the penalty term. We draw upon such penalty specifications from existing robust models, including Maenhout \cite{Maenhout2008}, Liu \cite{Liu2010}, and Pu and Zhang \cite{PZ2021}, among others. In particular, Liu \cite{Liu2010}, and Pu and Zhang \cite{PZ2021} focus explicitly on the Epstein-Zin utility function.
	\end{rem}

	To provide a concrete characterization of Assumption \ref{A3}, we further specify its implementation conditions as follows.
	
	\begin{ass}\label{A31}
		For the Epstein-Zin aggregator \eqref{nonunitEIS}, let the coefficients \( \gamma \) and \( \varphi \) satisfy one of the following four cases:
		
		\item[(i)] \( \gamma > 1 \) and \( \varphi > 1 \),
		
		\item[(ii)] \( \gamma > 1 \) and \( \varphi < 1 \), with \( \gamma \varphi \leq 1 \),
		
		\item[(iii)] \( \gamma < 1 \) and \( \varphi < 1 \),
		
		\item[(iv)] \( \gamma < 1 \) and \( \varphi > 1 \), with \( \gamma  \varphi \geq 1 \).
	\end{ass}

	As established in Proposition 3.2 of Kraft et al. \cite{KSS2013}, for the Epstein--Zin aggregator defined in \eqref{nonunitEIS}, Assumption \ref{A31} is sufficient for Assumption \ref{A3} to hold. In particular, in the unit EIS case as specified in \eqref{unitEIS}, Assumption \ref{A3} is clearly satisfied irrespective of whether \( \gamma > 1 \) or \( 0 < \gamma < 1 \).

	In view of the above analysis, we can directly apply Theorem \ref{VT} to derive the HJBI equation corresponding to problem (PR) \eqref{PR}.

	\begin{pro}\label{pro31}
		Let \( v \in C^{1,2}([0, T] \times (0, \infty) \times \mathbb{R}) \) be a solution of the HJBI equation
		\begin{equation}\label{HJBI1}
			\begin{aligned}
				&\max_{(\pi, q, c) \in \mathcal{U}_t(x, m)} \min_{\xi \in \mathcal{V}_t} \Bigg\{ f(c, v) + v_t + v_x[ rx + (\sigma m+a-r)\pi + \lambda \mu_1 \theta_1 q - c -  \sigma \xi_1 \pi  - \sqrt{\lambda \mu_2} \xi_3 q ] \\
				&\qquad \qquad \qquad \qquad - v_m(\alpha m + \beta \rho_1 \xi_1 + \beta \sqrt{1-\rho_1^2} \xi_2) +\frac{1}{2} v_{xx} (\sigma^2 \pi^2 + \lambda \mu_2 q^2) +\frac{1}{2} v_{mm} \beta^2 \\
				&\qquad \qquad \qquad \qquad + v_{xm} \sigma \beta \rho_1 \pi + \frac{1}{2\Psi} (\xi_1^2+\xi_2^2+\xi_3^2)\ \Bigg\} = 0, \\
				&v(T, x, m) =  \frac{1}{1 - \gamma}x^{1 - \gamma},
			\end{aligned}
		\end{equation}
		where the aggregator \( f \) is given by \eqref{nonunitEIS} or \eqref{unitEIS}. Assume
		
		\item[(i)] for any \( (\pi, q, c) \in \mathcal{U}_t(x, m) \) and \( \xi \in \mathcal{V}_t \), \eqref{RX},  \eqref{PR}   have a unique strong solution \( (X_s, J_s)_{t \leq s \leq T} \);
		
		\item[(ii)] there exist admissible \( \xi^* \in \mathcal{V}_t \) and \( (\pi^*, q^*, c^*) \in \mathcal{U}_t(x, m) \) satisfying for all \( (s, x, m) \in [t, T] \times \mathbb{R}^+ \times \mathbb{R} \),
		\begin{equation*}
			\begin{aligned}
				& f(c^*, v) + v_t + v_x[ rx + (\sigma m+a-r)\pi^* + \lambda \mu_1 \theta_1 q^* - c -  \sigma \xi_1^* \pi^*  - \sqrt{\lambda \mu_2} \xi_3^* q^* ] \\
				& - v_m(\alpha m + \beta \rho_1 \xi_1^* + \beta \sqrt{1-\rho_1^2} \xi_2^*) +\frac{1}{2} v_{xx} \left( \sigma^2 (\pi^*)^2 + \lambda \mu_2 (q^*)^2 \right) +\frac{1}{2} v_{mm} \beta^2 + v_{xm} \sigma \beta \rho_1 \pi^* \\
				& + \frac{1}{2\Psi} \left( (\xi_1^*)^2+(\xi_2^*)^2+(\xi_3^*)^2 \right) \ = 0,
			\end{aligned}
		\end{equation*}
		
		\item[(iii)] for the Epstein--Zin aggregator  defined in    \eqref{nonunitEIS}, Assumption \ref{A31} holds;
		
		\item[(iv)] the local martingale
		\begin{equation*}
			\int_t^\cdot v_x \sigma \pi(s) \textnormal{d}W_{1,s}^{\mathbb{Q}^\xi} + \int_t^\cdot v_x \sqrt{\lambda \mu_2} q(s) \textnormal{d}W_{3,s}^{\mathbb{Q}^\xi} + \int_t^\cdot v_m \beta\left( \rho_1 + \sqrt{1-\rho_1^2} \right) \textnormal{d}W_{2,s}^{\mathbb{Q}^\xi}
		\end{equation*}
		is a true martingale for any \( (\pi, q, c) \in \mathcal{U}_t(x, m) \) and \( \xi \in \mathcal{V}_t \);
		
		\item[(v)] for any \( (\pi, q, c) \in \mathcal{U}_t(x, m) \) and \( \xi \in \mathcal{V}_t \), we have
		\begin{equation*}
			\mathbb{E}\left[ \int_t^T |v(s, X_s) - J_s| \, \textnormal{d}s \right] < \infty.
		\end{equation*}
		Then \( (\pi^*, q^*, c^*, \xi^*) \) is an optimal control and \( v \) is the value function of problem (PR) \eqref{PR}.
		
	\end{pro}

	To obtain an explicit solution to the model, we adopt the methodology proposed by Maenhout \cite{Maenhout2004, Maenhout2008}. In particular, we assume that the robustness preference parameter \( \Psi \) is state-dependent and is scaled by the value function as follows:
	\begin{equation}\label{robpar}
		\Psi(t, x, m) = \frac{\Phi}{(1-\gamma) v(t, x, m)},
	\end{equation}
	where \( \Phi \geq 0 \) is a constant. This parameter \( \Phi \) represents the ambiguity-aversion coefficient describing the insurer's attitudes to diffusion.

	\section{The robust optimal strategy}\label{sec4}

	\subsection{Optimal strategy under EIS $\neq$ 1}\label{sec4.1}

	In this subsection, for the Epstein--Zin aggregator \( f \)  defined in \eqref{nonunitEIS}, we derive the optimal robust reinsurance--investment--consumption strategy via the first-order condition (FOC) method and propose a candidate solution \( v(t, x, m) \) for the HJBI equation \eqref{HJBI1}. From \eqref{HJBI1}, the FOC with respect to \( \xi \) yields
	\begin{equation}\label{xi*}
		\xi^* = \left(\xi_1^*, \xi_2^*, \xi_3^* \right) = \left(\frac{(v_x \pi \sigma + v_m \beta \rho_1)\Phi}{(1-\gamma)v}, \frac{v_m \beta \sqrt{1-\rho_1^2} \Phi}{(1-\gamma)v}, \frac{v_x \sqrt{\lambda \mu_2} q \Phi}{(1-\gamma)v}\right).
	\end{equation}
	
	Substituting \eqref{xi*} back into the HJBI equation \eqref{HJBI1}, we obtain
	\begin{equation}\label{HJBI2}
		\begin{aligned}
			&\max_{(\pi, q, c) \in \mathcal{U}_t(x, m)}  \Bigg\{ f(c, v) + v_t + v_x \Big[ rx + (\sigma m+a-r)\pi + \lambda \mu_1 \theta_1 q - c - \frac{\Phi v_x \sigma^2 \pi^2}{(1-\gamma)v} \\
			&\qquad \qquad \qquad -\frac{2 \Phi v_m \beta \rho_1 \sigma \pi}{(1-\gamma)v} - \frac{\Phi v_x \lambda \mu_2 q^2}{(1-\gamma)v} \Big] - \Big[ \alpha m v_m + \frac{\Phi v_m^2 \beta^2}{(1-\gamma)v} \Big] \\
			&\qquad \qquad \qquad  + \frac{1}{2} v_{xx} (\sigma^2 \pi^2 + \lambda \mu_2 q^2) +\frac{1}{2} v_{mm} \beta^2 + v_{xm} \sigma \beta \rho_1 \pi \\
			&\qquad \qquad \qquad  + \frac{\Phi}{2(1-\gamma)v} (v_x^2 \sigma^2 \pi^2 + 2v_x v_m \sigma \beta \rho_1 \pi + v_m^2 \beta^2 + v_x^2 \lambda \mu_2 q^2)\ \Bigg\} = 0
		\end{aligned}
	\end{equation}
	with the terminal condition \( v(T, x, m) = \frac{1}{1-\gamma}x^{1-\gamma} \).
	Noting that the terms involving \( \pi, q \) and \( c \) in equation \eqref{HJBI2} are separable, we decompose these components and apply the FOC to obtain
	\begin{equation}\label{piqc*}
		\begin{aligned}
			\pi^* & = \frac{v_x (\sigma m+a-r) + v_{xm} \beta \rho_1 \sigma - \frac{v_x v_m \Phi \beta \rho_1 \sigma}{(1-\gamma)v}}{\left( \frac{v_x^2 \Phi}{(1-\gamma)v} - v_{xx} \right) \sigma^2}, \\
			q^* & = \frac{v_x \theta_1 \mu_1}{\frac{v_x^2 \mu_2 \Phi}{(1-\gamma)v} - v_{xx} \mu_2}, \\
			c^* & = v_x^{-\varphi} \delta^\varphi \left[(1-\gamma) v \right]^{\frac{1-\gamma \varphi}{1-\gamma}}.
		\end{aligned}
	\end{equation}

	Utilizing a constant \( k \), we construct the conjectured form as:
	\begin{equation}\label{conject}
		v(t, x, m) = \frac{1}{1-\gamma}x^{1-\gamma} g(t, m)^k, \quad (t, x, m) \in [0, T] \times \mathbb{R}^+ \times \mathbb{R},
	\end{equation}
	with \( g \in C^{1,2}([0, T] \times \mathbb{R}) \) and \( g(T, m) = 1 \).

	In order to reduce the nonlinear partial differential equation (PDE) \eqref{HJBI2} to a linear PDE under the conjectures \eqref{robpar} and \eqref{conject}, we impose the following assumptions on the coefficients \( k \) and \( \varphi \) as in Pu and Zhang \cite{PZ2021}. These assumptions can also be derived from the calculations in Proposition \ref{pro-exact}.
	
	\begin{ass}\label{kvarphi}
		The coefficients \( k \) and \( \varphi \) satisfy
		\begin{equation*}
			k = 1 \Big/ \left( 1 - \frac{\Phi}{1 - \gamma} + \frac{(1 - \gamma - \Phi)^2}{(1 - \gamma)(\Phi + \gamma)} \rho_1^2  \right)
		\end{equation*}
		and
		\begin{equation*}
			\varphi = 2 - \gamma - \Phi + \frac{(1 - \gamma - \Phi)^2}{\Phi + \gamma} \rho_1^2.
		\end{equation*}
	\end{ass}

	The following result provides a solution to Problem (PR) \eqref{PR}.
	
	\begin{pro}\label{pro-exact}
		Under Assumption \ref{kvarphi} and conditions (i)--(ii) of Proposition \ref{pro31}, the candidate solution to  the HJB equation \eqref{HJBI2} is given by
		\begin{equation}\label{conject1}
			v(t, x, m) = \frac{1}{1 - \gamma} x^{1 - \gamma} g(t, m)^k, \qquad (t, x, m) \in [0, T] \times \mathbb{R}^+ \times \mathbb{R},
		\end{equation}
		where the function \( g \) satisfies the following PDE:
		\begin{equation}\label{g1}
			\begin{aligned}
				& g_t + \frac{1-\gamma}{k} \left[ r + \frac{(\sigma m + a - r)^2}{2(\Phi+\gamma)\sigma^2} + \frac{\lambda\theta_1^2\mu_1^2}{2(\Phi+\gamma)\mu_2} - \frac{\delta}{1-\frac{1}{\varphi}} \right] g  \\
				& + \left[ \frac{(1-\gamma-\Phi)}{(\Phi+\gamma)\sigma}(\sigma m + a - r)\beta \rho_1 - \alpha m \right] g_m + \frac{1}{2}\beta^2 g_{mm} = 0.
			\end{aligned}
		\end{equation}
		Moreover, for this \(v(t, x, m)\), if we further assume the conditions (iv)--(v) of Proposition \ref{pro31}, the optimal reinsurance-investment-consumption strategy is
		\begin{equation}\label{piqc**}
			\begin{aligned}
				\pi^* & = \frac{x}{(\Phi+\gamma)\sigma^2} \left[(\sigma m + a - r)  + \left(1-\frac{\Phi}{1-\gamma}\right) k \beta \rho_1 \sigma \frac{g_m}{g} \right], \\
				q^* & = x \frac{\theta_1 \mu_1}{(\Phi+\gamma) \mu_2}, \\
				c^* & = x \delta^\varphi g^{ k \frac{1-\varphi}{1-\gamma}}.
			\end{aligned}
		\end{equation}
	\end{pro}

	For \(g\) in the value function \eqref{conject1}, we characterize it by means of the Feynman-Kac representation theorem, which is analogous to Theorem 4.1 in Kraft et al. \cite{KSS2013}.
	Define
	\begin{equation*}
		\begin{aligned}
			& H_1 := \frac{1-\gamma}{k} \left(r + \frac{(\sigma m+a-r)^2}{2(\Phi+\gamma)\sigma^2} + \frac{\lambda \theta_1^2 \mu_1^2}{2(\Phi+\gamma)\mu_2} - \frac{\delta}{1-\frac{1}{\varphi}}\right), \\
			& H_2 := \frac{(1-\gamma-\Phi)}{(\Phi+\gamma)\sigma} (\sigma m+a-r) \beta \rho_1 - \alpha m.
		\end{aligned}
	\end{equation*}
	The corresponding results are summarized in the following proposition.
	
	\begin{pro}\label{proFK}
		If the coefficients in \eqref{g1} satisfy the conditions of Heath and Schweizer \cite{Heath2000}, then the solution \(g\) admits the Feynman-Kac representation
		\begin{equation*}
			g(t, m) = \delta^\varphi H(t, m) +  h(t, m; T),
		\end{equation*}
		where
		\begin{equation*}
			H(t, m) := \int_t^T h(t, m; s) \textnormal{d}s, \
			h(t, m; s) := \tilde{\mathbb{E}}_{t,m}\left[ e^{ \int_t^s H_1(m(u)) \textnormal{d}u } \right]
		\end{equation*}
		for \(t \in [0, T]\), \(m \in \mathbb{R}\), and the expectation \(\tilde{\mathbb{E}}_{t,m}[\cdot]\) is taken with respect to the equivalent measure \(\tilde{\mathbb{P}}\), under which the process \(m\) is conditioned on \(m(t) = m\) and has drift \(H_2 \). Here we also note that \(h\) satisfies the PDE
		\begin{equation}\label{h1}
			h_t + H_1 h + H_2 h_m + \frac{1}{2}\beta^2 h_{mm} = 0 \quad \text{on } [0, s] \times \mathbb{R}
		\end{equation}
		with the terminal condition \(h(s, m; s) = 1\).
	\end{pro}

	Beyond presenting the candidate optimal strategy, we comprehensively validate that conditions (iv) and (v) required in Proposition \ref{pro31} are satisfied. Furthermore, in accordance with the definitions of \(\mathcal{U}_t(x, m)\) and \(\mathcal{V}_t\) provided in Section \ref{sec3}, we introduce the following set of admissible controls: 
	\begin{equation*}
		\begin{aligned}
			\mathcal{A}_t := \Big\{ &(\pi, q, c;\xi): (\pi,q, c) \in \mathcal{U}_t(x, m), \, \xi \in \mathcal{V}_t, \, \frac{c}{x}(s, m) \leq \frac{1}{T - s} + b_0 m ^2+ b_1 |m| + b_2; \\
			& \pi = \tilde{\pi} x, \, |\tilde{\pi}| \leq D_1 |m| + D_2; \, q = \tilde{q} x, \, \tilde{q} \, \text{is a constant}, \, 0 \leq \tilde{q} \leq D_3;  \\
			& \xi = (\xi_1, \xi_2, \xi_3), \, |\xi_1| \leq D_4 |m| + D_5, \, \xi_2 \, \text{and }\, \xi_3 \, \text{have an upper bound with } 0 \leq \xi_3 \leq D_6; \\
			& X \, \text{satisfies  condition (M)} , \, X_s \geq 0, \, s \in [t, T] \Big\},
		\end{aligned}
	\end{equation*}
	where $b_0, \, b_1, \, b_2, \, D_i, \, 1 \leq i \leq 6$ are constants, and the condition (M) is
	\begin{equation*}
		\mathbb{E}^{\mathbb{Q}^\xi}[X_t^{-\ell}] \leq K (T - t)^{-\ell}, \quad t \in [0, T], \text{ for some } K > 0, \, \ell > 2(\gamma - 1). \tag{M}
	\end{equation*}
	Moreover, we impose the following additional assumptions on the coefficients:
	\begin{itemize}
		\item[(H1)] \(1 < \gamma < \min\left\{k+\frac{3}{2}, \frac{1}{\overline{k}} + 1\right\}\) where \(\overline{k} > 2\) such that \(\overline{k}-2\) is sufficiently small, and \(\rho_1 \leq 0\);
		\item[(H2)] \( \alpha (\Phi+\gamma)^2 > 8 \beta^2 T \cdot \max\left\{\Phi^2, \, \overline{k}^2 (\gamma-1)^2 \right\} \);
		\item[(H3)] \( \alpha > 16 \overline{k} (\gamma-1) \beta^2 T \left( \frac{ 2k + (\gamma-1)\sigma }{ k(\Phi+\gamma) \sigma } + \frac{2\Phi + \overline{k} (\gamma-1) + 1}{(\Phi+\gamma)^2} \right) + \frac{(\overline{k}(1-\gamma) - \Phi) \beta \rho_1}{\Phi+\gamma}\).
	\end{itemize}

	\begin{thm}\label{thm-exact}
		Assume that  conditions  (i)--(ii) of Proposition \ref{pro31}, Assumptions \ref{A31}, \ref{kvarphi} and (H1)--(H3)
		hold. Then, the  value function $v$ in our model admits the
		representation \eqref{conject}, where
		\begin{equation}\label{g*}
			g(t, m) = \delta^\varphi \int_t^T e^{A(t,s)-B(t,s)m-C(t,s)m^2} \, \textnormal{d}s + e^{A(t,T)-B(t,T)m-C(t,T)m^2},
		\end{equation}
		and the functions \( A \), \( B \) and \( C \) are given by
		\begin{equation}\label{ABCsolution}
			\begin{aligned}
				A(t, s) &= \int_{t}^{s} \left[ \frac{1}{2} \beta^2 B^2(\tau,s) - \beta^2 C(\tau,s) - \frac{1-\gamma-\Phi}{(\Phi+\gamma)\sigma} (a-r) \beta \rho_1 B(\tau,s)  \right] \textnormal{d}\tau  \\
				&\quad + \frac{1-\gamma}{k}  \left[ r + \frac{(a-r)^2}{2(\Phi+\gamma)\sigma^2} + \frac{\lambda \theta_1^2 \mu_1^2}{2(\Phi+\gamma)\mu_2} - \frac{\delta}{1-\frac{1}{\varphi}} \right] (s-t), \\
				B(t, s) &= \frac{(1-\gamma)(a-r)}{k(\Phi+\gamma)\sigma} \int_t^s \left( 2\frac{1-\gamma-\Phi}{1-\gamma} k \beta \rho_1 C(\tau, s) - 1\right)  e^{-\int_t^\tau \left( \kappa + 2\beta^2 C(u, s) \right) \textnormal{d}u} \, \textnormal{d}\tau, \\
				C(t, s) &= 2b_0 \frac{e^{\Delta(s-t)} - 1}{e^{\Delta(s-t)}(2\kappa + \Delta) - 2\kappa + \Delta},
			\end{aligned}
		\end{equation}
		with \( \kappa := \alpha - \frac{1-\gamma-\Phi}{\Phi+\gamma} \beta \rho_1 \), \( b_0 := -\frac{1-\gamma}{2k(\Phi+\gamma)} \) and \( \Delta := 2 \sqrt{\kappa^2 + 2\beta^2 b_0} \). The optimal reinsurance-investment-consumption strategy is
		\begin{equation}\label{piqc***}
			\begin{aligned}
				\pi^*(t) & = \frac{x(t)}{(\Phi+\gamma)\sigma^2} \left[(\sigma m(t)+a-r)  + \frac{1-\gamma-\Phi}{1-\gamma} k \beta \rho_1 \sigma \frac{g_m(t,m(t))}{g(t,m(t))} \right], \\
				q^*(t) & = x(t) \frac{\theta_1 \mu_1}{(\Phi+\gamma) \mu_2}, \\
				c^*(t) & = x(t) \delta^\varphi g^{-1}(t,m(t)),
			\end{aligned}
		\end{equation}
		and the distortion process \( \xi^* = \left(\xi_1^*, \xi_2^*, \xi_3^* \right)\) is
		\begin{equation}\label{xi**}
			\begin{aligned}
				\xi_1^*(t) & = \frac{\Phi}{(\Phi+\gamma)\sigma} (\sigma m(t)+a-r) + \frac{\Phi k \beta \rho_1}{(\Phi+\gamma) (1-\gamma)} \frac{g_m(t,m(t))}{g(t,m(t))}, \\
				\xi_2^*(t) & = \frac{\Phi k \beta \sqrt{1-\rho_1^2}}{1-\gamma}  \frac{g_m(t,m(t))}{g(t,m(t))}, \\
				\xi_3^*(t) & = \frac{\Phi \theta_1 \mu_1 \sqrt{\lambda}}{(\Phi+\gamma) \sqrt{\mu_2}}.
			\end{aligned}
		\end{equation}
	\end{thm}

	\begin{rem}
		Under appropriate conditions, our results  reduce to those of Proposition 5.3 in Kraft et al. \cite{KSS2013}. Indeed, setting $\Phi =0$ and $a = r$ in \eqref{piqc***} makes the resulting  $\tilde{\pi}^*$ and $\tilde{c}^*$  identical to $\pi_t^*$ and $\frac{c_t^*}{X_t^*}$ in Kraft et al. \cite[Proposition 5.3]{KSS2013}, highlighting the generality and inclusiveness of our framework. Furthermore, a key innovation of this study is the rigorous proof that our constructed strategy satisfies all conditions of the verification theorem.
	\end{rem}

	\subsection{Optimal strategy under EIS = 1}
	
	In this subsection, we present an exact solution for the Epstein--Zin aggregator \eqref{unitEIS} with EIS \( = 1 \).
	Analogously to Subsection \ref{sec4.1}, we derive \eqref{xi*} and the corresponding HJB equation \eqref{HJBI2} in the same form, with
	\( f \) given by \eqref{unitEIS}.
	Consequently, we obtain \( \pi^* \), \( q^* \), and \( c^* \) in the following \underline{}form
	\begin{equation*}\label{piqc*unit}
		\begin{aligned}
			\pi^* & = \frac{v_x (\sigma m + a - r) + v_{xm} \beta \rho_1 \sigma - \frac{v_x v_m \Phi \beta \rho_1 \sigma}{(1-\gamma)v}}{\left( \frac{v_x^2 \Phi}{(1-\gamma)v} - v_{xx} \right) \sigma^2}, \\
			q^* & = \frac{v_x \theta_1 \mu_1}{\frac{v_x^2 \mu_2 \Phi}{(1-\gamma)v} - v_{xx} \mu_2}, \\
			c^* & = \delta (1-\gamma) \frac{v}{v_x}.
		\end{aligned}
	\end{equation*}

	We construct the conjectured form as:
	\begin{equation*}\label{conjectEIS}
		v(t, x, m) = \frac{1}{1-\gamma}x^{1-\gamma} g(t, m), \quad (t, x, m) \in [0, T] \times \mathbb{R}^+ \times \mathbb{R},
	\end{equation*}
	with \( g \in C^{1,2}([0, T] \times \mathbb{R}) \) and \( g(T, m) = 1 \).
	The following result provides a solution to Problem (PR) \eqref{PR}.
	
	\begin{pro}
		If conditions $(i)$, $(ii)$, $(iv)$ and $(v)$ of Proposition~\ref{pro31} are satisfied, the candidate solution to the HJBI equation~\eqref{HJBI1} is given by
		\begin{equation*}\label{conject1EIS}
			v(t, x, m) = \frac{1}{1 - \gamma} x^{1 - \gamma} g(t, m), \qquad (t, x, m) \in [0, T] \times \mathbb{R}^+ \times \mathbb{R}.
		\end{equation*}
		Here, the function \( g \) is
		\begin{equation}\label{conjectg}
			g(t, m) =  e^{G(t)m^2 + L(t)m + H(t)},
		\end{equation}
		where the functions \( G(t) \), \( L(t) \) and \( H(t) \) are given by
		\begin{equation*}\label{GLHsolution}
			\begin{aligned}
				G(t) &= 2G_3 \frac{e^{\sqrt{G_2^2-4G_1G_3}(T-t)} - 1}{e^{\sqrt{G_2^2-4G_1G_3}(T-t)}(G_2 + \sqrt{G_2^2-4G_1G_3}) - G_2 + \sqrt{G_2^2-4G_1G_3}}, \\
				L(t) &= \frac{(1-\gamma)(a-r)}{(\Phi+\gamma)\sigma} \int_{t}^{T} \left( 2 \frac{1-\gamma-\Phi}{1-\gamma} \beta \rho_1 G(s) + 1 \right)  e^{\int_{t}^{s} \left( \frac{1-\gamma-\Phi}{\Phi+\gamma} \beta \rho_1 -\alpha -\delta - G_1 G(u) \right) \textnormal{d}u} \, \textnormal{d}s, \\
				H(t) &= \int_{t}^{T} \left[ \frac{1}{2} (\beta^2 + 2G_0) L^2(s) + \frac{(1-\gamma-\Phi)(a-r)}{(\Phi+\gamma)\sigma} \beta \rho_1 L(s) + \beta^2 G(s) \right] e^{\delta (s-t)} \textnormal{d}s  \\
				&\quad + \frac{1-\gamma}{\delta}  \left[ r + \frac{(a-r)^2}{2(\Phi+\gamma)\sigma^2} + \frac{\lambda \theta_1^2 \mu_1^2}{2(\Phi+\gamma)\mu_2} + \delta (\ln \delta - 1) \right] \left( e^{\delta(T-t)} - 1 \right),
			\end{aligned}
		\end{equation*}
		with \( G_1 := -2(\beta^2 + 2G_0) \), \( G_2 := 2\alpha + \delta - \frac{2(1-\gamma-\Phi)\beta \rho_1}{\Phi+\gamma} \), \( G_3 := - \frac{1-\gamma}{2(\Phi+\gamma)} \) and \( G_0 := \frac{(1-\gamma-\Phi)^2 \beta^2 \rho_1^2}{2 (\Phi+\gamma) (1-\gamma)^2} - \frac{\Phi \beta^2}{2(1-\gamma)} \). The optimal  strategies  are
		\begin{equation*}\label{piqc***EIS}
			\begin{aligned}
				\pi^*(t) & = \frac{x(t)}{(\Phi+\gamma)\sigma^2} \left[(\sigma m(t) + a - r) +
				\frac{1-\gamma-\Phi}{1-\gamma} k \beta \rho_1 \sigma \left(2 G(t) m(t) + L(t)\right) \right], \\
				q^*(t) & = \frac{\theta_1 \mu_1}{(\Phi+\gamma) \mu_2} x(t) , \\
				c^*(t) & = \delta x(t),
			\end{aligned}
		\end{equation*}
		and the distortion process \( \xi^* = \left(\xi_1^*, \xi_2^*, \xi_3^* \right)\) is
		\begin{equation*}\label{xi*EIS}
			\begin{aligned}
				\xi_1^*(t) & = \frac{\Phi}{(\Phi+\gamma)\sigma} (\sigma m(t) + a - r) + \frac{\Phi k \beta \rho_1}{(\Phi+\gamma) (1-\gamma)} \left(2 G(t) m(t) + L(t)\right), \\
				\xi_2^*(t) & = \frac{\Phi k \beta \sqrt{1-\rho_1^2}}{1-\gamma}  \left(2 G(t) m(t) + L(t)\right), \\
				\xi_3^*(t) & = \frac{\Phi \theta_1 \mu_1 \sqrt{\lambda}}{(\Phi+\gamma) \sqrt{\mu_2}}.
			\end{aligned}
		\end{equation*}
		
	\end{pro}

	\subsection{Numerical analysis}

	In this subsection, we will explore the potential financial and insurance implications of the findings in this study and conduct a parameter sensitivity analysis of these strategies. Unless otherwise specified, the model parameters are as shown in Table \ref{para}, with some values drawn from the empirical study of Liu and Pan \cite{LP2003} and the assumptions of Kraft et al. \cite{KSS2013}.
	Additionally, we assume that the claim size distribution follows a Gamma distribution, specifically $\Gamma(4, 0.25)$.

	\renewcommand{\arraystretch}{2}
	\begin{table}[H]
		\scriptsize
		\centering
		\caption{Values of model parameters}
		\label{para}
		\begin{tabular}{l l l}
			\hline
			\textbf{Epstein--Zin preference parameters} & \textbf{Ambiguity-aversion coefficient} & \textbf{Reinsurance safety loading} \\ \hline
			$\gamma = 1.2$, $\delta = 0.08$ & $\Phi = 0.8$ & $\theta_1 = 0.2$ \\ \hline
			\textbf{Risky asset parameters} & \textbf{Risk-free asset parameter} & \textbf{} \\ \hline
			$\sigma = 0.2$, $\beta = 0.25$, $\rho_1 = -0.5$, $\alpha = 5$, $a = 0.07$ & $r = 0.02$ &  \\ \hline
			\textbf{Claim amount parameters} & \textbf{Claim-number intensity} & \textbf{Auxiliary parameter} \\ \hline
			$\mu_1 = 1.0$, $\mu_2 = 1.25$ & $\lambda = 1.0$ & $\bar{k} = 2.1$ \\ \hline
			\textbf{Initial time} & \textbf{Terminal time} & \textbf{} \\ \hline
			$t = 0.5$ & $T = 1.0$ &  \\ \hline
		\end{tabular}
	\end{table}

\begin{figure}[h]
		\begin{minipage}[b]{0.53\linewidth}
			\centering
			\includegraphics[width=\textwidth]{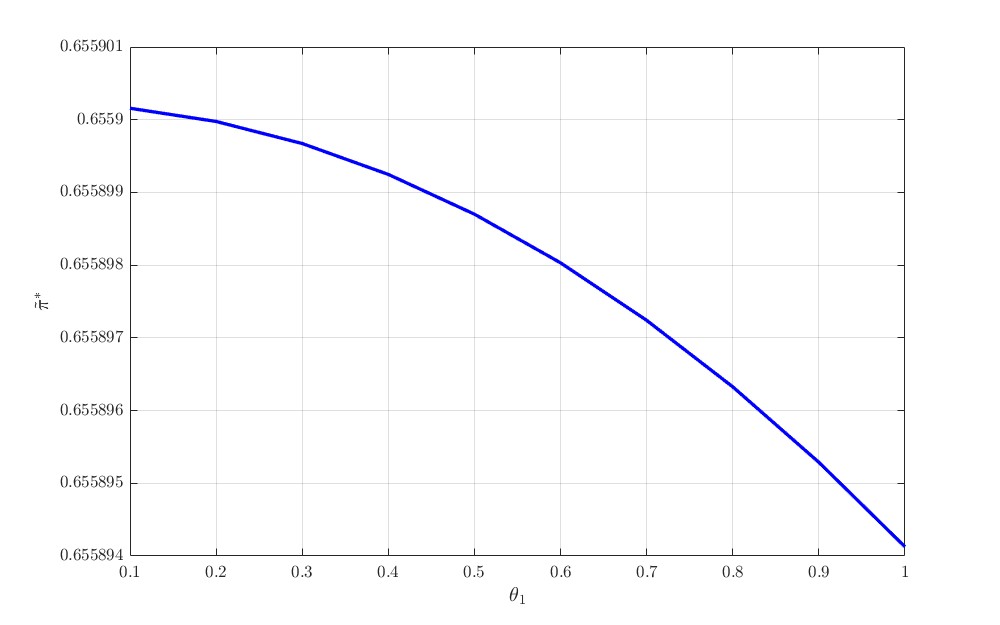}
			\centerline{\small(a) the investment--wealth ratio $\tilde{\pi}^*$}
		\end{minipage}
		\hfill
		\begin{minipage}[b]{0.53\linewidth}
			\centering
			\includegraphics[width=\textwidth]{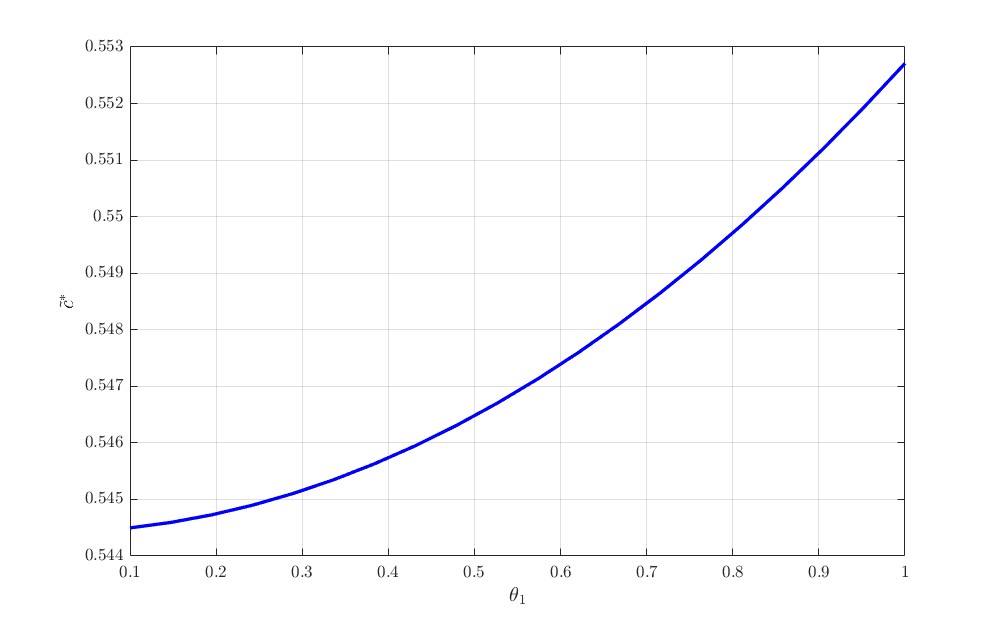}
			\centerline{\small(b) the consumption--wealth ratio $\tilde{c}^*$}
		\end{minipage}
		
		\vspace{1mm}
		\caption{The effect of the reinsurance safety loading coefficient $\theta_1$}
		\label{fig1}
	\end{figure}

	How the reinsurance safety loading coefficient affects the insurer's optimal investment strategy and internal surplus allocation (consumption), given its reinsurance business, is a key question worth exploring. This is illustrated in Figures \ref{fig1} (a) and (b), which show the investment--wealth ratio $\tilde{\pi}^*$ and the consumption--wealth ratio $\tilde{c}^*$, respectively, as functions of
	$\theta_1$ over the range of 0.1 to 1.0.

	\setlength{\abovecaptionskip}{2pt}
	\vspace{1mm}

	It can be observed that $\tilde{\pi}^*$ decreases monotonically with $\theta_1$. The underlying mechanism is that as $\theta_1$ increases, the insurance company tends to reduce its reinsurance coverage, i.e., to retain a larger share of the risk, thereby increasing its risk exposure. To hedge against this additional risk, the insurer correspondingly reduces the amount invested in the risky asset. Notably, $\tilde{\pi}^*$ is not highly sensitive to changes in $\theta_1$, exhibiting only a mild downward trend. This is mainly because variations in reinsurance pricing have a limited impact on the attractiveness of investing in the risky asset. In contrast, $\tilde{c}^*$ increases monotonically with $\theta_1$. This is primarily driven by the reduction in total reinsurance expenditure. The resulting savings in reinsurance premiums are directly converted into current profits, releasing more distributable surplus and consequently raising the level of surplus allocation.

	In the following, we compare the optimal investment--wealth ratios for robust insurers $(\Phi=0.4, \Phi=0.8)$ and a non-robust insurer $(\Phi=0)$.

	\setlength{\abovecaptionskip}{2pt}
	\begin{figure}[h]
		\centering
		\subfigure{
			\includegraphics[scale=0.3]{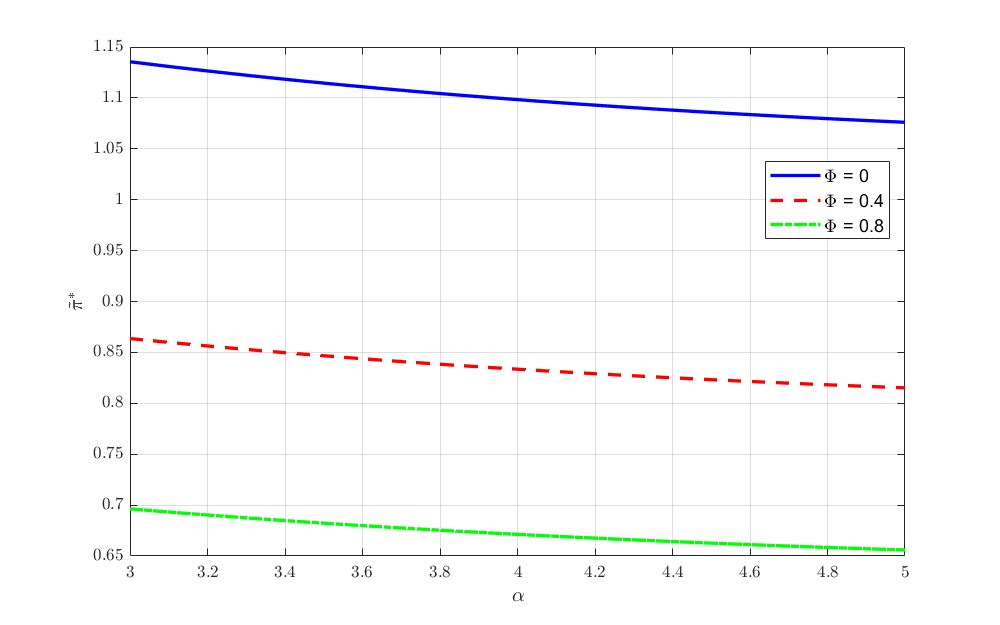}}
		\caption{The effect of the mean reversion rate $\alpha$ on the investment--wealth ratio $\tilde{\pi}^*$}
		\label{fig3}
	\end{figure}
	\vspace{1mm}

	Figure \ref{fig3} shows that the optimal investment--wealth ratio $\tilde{\pi}^*$ decreases with $\alpha$ under both non-robustness $(\Phi=0)$ and robustness $(\Phi=0.4, \Phi=0.8)$, consistent with real market behavior. When $\alpha$ is small, strong drift persistence generates stable long-term excess return expectations, leading to a higher investment–wealth ratio. When $\alpha$ increases, drift persistence weakens, reducing the predictive power for future returns and shortening the window for capturing excess returns, thereby lowering the ratio.

	\setlength{\abovecaptionskip}{2pt}
	\begin{figure}[h]
		\centering
		\subfigure{
			\includegraphics[scale=0.3]{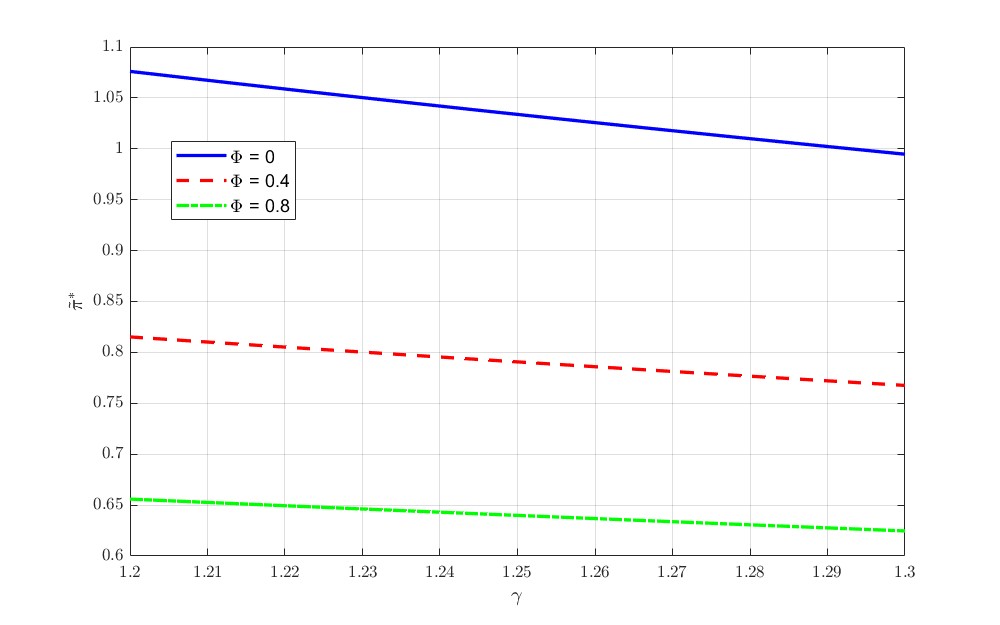}}
		\caption{The effect of the risk aversion coefficient $\gamma$ on the investment--wealth ratio $\tilde{\pi}^*$}
		\label{fig4}
	\end{figure}
	\vspace{1mm}

	Figure \ref{fig4} demonstrates that $\tilde{\pi}^*$ is negatively correlated with the risk aversion coefficient $\gamma$, which is consistent with the fundamental principles of insurance financial markets. Specifically, an increase in $\gamma$ strengthens the insurer's risk-averse behavior, leading to a reduction in the allocation to risky assets.

	Figures \ref{fig3} and \ref{fig4} reveal that, for any fixed $\alpha$  or $\gamma$, more robust insurers allocate less to the risky asset than non-robust insurers. This is consistent with the theory that a larger ambiguity-aversion coefficient reflects greater insurer concern about model misspecification, thereby reducing risky asset allocation. Moreover, as $\Phi$ increases from 0.4 to 0.8, the investment--wealth ratio curve shifts further downward, indicating that a stronger preference for robustness further dampens risk-taking.

	\section{Campbell--Shiller approximation under EIS $\neq$ 1}\label{sec5}

	The Campbell--Shiller approximation is a landmark tool in financial economics, whose core value lies in its ability to transform complex nonlinear equations into a tractable linear framework, thereby facilitating theoretical derivations.
	In this section, we employ this approximation method to analyze the problem under EIS $\neq$ 1, and compare the obtained results with the exact solution derived in Subsection \ref{sec4.1}.

	Following the equation \eqref{HJBI3} in Appendix and multiplying both sides of this equation by \( g(1-\gamma) k^{-1} \), we obtain
	\begin{equation}\label{HJBI5}
		\begin{aligned}
			& g_t + \frac{1-\gamma}{k} \left[\left(1 - \frac{1}{\varphi}\right)^{-1} - 1\right] \delta^{\varphi} g^{k\frac{1-\varphi}{1-\gamma} + 1} + \frac{1-\gamma}{k}  \Bigg[ \frac{(\sigma m + a - r)^2}{2(\Phi+\gamma)\sigma^2} + \frac{\lambda \theta_1^2 \mu_1^2}{2(\Phi+\gamma)\mu_2} \\
			& + r - \frac{\delta}{1-\frac{1}{\varphi}} \Bigg] g + \frac{1}{2}\beta^2 g_{mm} + \left[\frac{1-\gamma-\Phi}{(\Phi+\gamma)\sigma} (\sigma m + a - r) \beta \rho_1 - \alpha m \right] g_m \\
			& + \left[\frac{(1-\gamma-\Phi)^2}{2(\Phi+\gamma)(1-\gamma)} k \beta^2 \rho_1^2 - \frac{\Phi k \beta^2}{2(1-\gamma)} + \frac{\beta^2 (k-1)}{2} \right] \frac{g_m^2}{g} = 0.
		\end{aligned}
	\end{equation}
	
	Applying a first-order log-linear expansion of the consumption--wealth ratio around its long-run steady state, i.e.,
	\begin{equation}\label{CS}
		\frac{c}{x} = e^{\ln \left(\frac{c}{x}\right)} \approx w \left( 1 - \ln w + \ln \frac{c}{x} \right) = w \left( 1 - \ln w + \ln \delta^{\varphi} + k \frac{1-\varphi}{1-\gamma} \ln g \right),
	\end{equation}
	where
	\(
	\ln w = \mathbb{E}^{\mathbb{Q}^\xi} \left[ \ln \left( \frac{c}{x} \right) (m_\infty) \right]
	\)
	and \( m_\infty \) is a random variable that has the stationary distribution of the process \( m \).
	Following Kraft et al. \cite{KSS2013} and Han and Hung \cite{HH2017}, we treat \( w \) as a fixed constant. Substituting the approximation \eqref{CS} into the HJB equation \eqref{HJBI5} yields the differential equation for the Campbell--Shiller approximation \( g^{\text{CS}} \) of \( g \):
	\begin{equation}\label{HJBI6}
		\begin{aligned}
			& g_t^{\text{CS}} + \frac{1-\gamma}{k(\varphi-1)} w \left( 1 - \ln w + \varphi \ln \delta + k \frac{1-\varphi}{1-\gamma} \ln g \right) g^{\text{CS}} + \frac{1-\gamma}{k}  \Bigg[ \frac{(\sigma m + a - r)^2}{2(\Phi+\gamma)\sigma^2}  \\
			& + \frac{\lambda \theta_1^2 \mu_1^2}{2(\Phi+\gamma)\mu_2} + r - \frac{\delta}{1-\frac{1}{\varphi}} \Bigg] g^{\text{CS}} + \frac{1}{2}\beta^2 g_{mm}^{\text{CS}} + \left[\frac{1-\gamma-\Phi}{(\Phi+\gamma)\sigma} (\sigma m + a - r) \beta \rho_1 - \alpha m \right] g_m^{\text{CS}} \\
			& + \left[\frac{(1-\gamma-\Phi)^2}{2(\Phi+\gamma)(1-\gamma)} k \beta^2 \rho_1^2 - \frac{\Phi k \beta^2}{2(1-\gamma)} + \frac{\beta^2 (k-1)}{2} \right] \frac{(g_m^{\text{CS}})^2}{g^{\text{CS}}} = 0.
		\end{aligned}
	\end{equation}
	In this case, letting \( g^{\text{CS}} \) take the form as in \eqref{conjectg}, we can  obtain the solution to equation \eqref{HJBI6}. The results are summarized in the following theorem.

	\begin{thm}
		Under  the Campbell--Shiller approximation approach, the candidate solution to the HJBI equation~\eqref{HJBI1} is given by
		\begin{equation*}\label{conject1CS}
			v^{\textnormal{CS}}(t, x, m) = \frac{1}{1 - \gamma} x^{1 - \gamma} \left(g^{\textnormal{CS}}(t, m)\right)^k, \qquad (t, x, m) \in [0, T] \times \mathbb{R}^+ \times \mathbb{R}.
		\end{equation*}
		Here, the function \( g^{\textnormal{CS}} \) is
		\begin{equation*}\label{conjectgCS}
			g^{\textnormal{CS}}(t, m) =  e^{G(t)m^2 + L(t)m + H(t)},
		\end{equation*}
		where the functions \( G(t) \), \( L(t) \) and \( H(t) \) are given by
		\begin{equation*}\label{GLHsolutionCS}
			\begin{aligned}
				G(t) &= 2G_3 \frac{e^{\sqrt{G_2^2-4G_1G_3}(T-t)} - 1}{e^{\sqrt{G_2^2-4G_1G_3}(T-t)}(G_2 + \sqrt{G_2^2-4G_1G_3}) - G_2 + \sqrt{G_2^2-4G_1G_3}}, \\
				L(t) &= \frac{(1-\gamma)(a-r)}{k(\Phi+\gamma)\sigma} \int_{t}^{T} \left( 2 \frac{1-\gamma-\Phi}{1-\gamma} k \beta \rho_1 G(s) + 1 \right)  e^{\int_{t}^{s} \left( \frac{1-\gamma-\Phi}{\Phi+\gamma} \beta \rho_1 - \alpha - w - G_1 G(u) \right) \textnormal{d}u} \, \textnormal{d}s, \\
				H(t) &= \int_{t}^{T} \left[ \frac{1}{2} (\beta^2 + 2G_0) L^2(s) + \frac{(1-\gamma-\Phi)(a-r)}{(\Phi+\gamma)\sigma} \beta \rho_1 L(s) + \beta^2 G(s) \right] e^{\delta (s-t)} \textnormal{d}s  \\
				&\quad + \frac{1-\gamma}{k} \left[ \frac{w}{\varphi-1} \left(1- \ln w + \varphi \ln \delta\right) + r - \frac{\delta}{1-\frac{1}{\varphi}} + \frac{(a-r)^2}{2(\Phi+\gamma)\sigma^2} + \frac{\lambda \theta_1^2 \mu_1^2}{2(\Phi+\gamma)\mu_2} \right] \left( e^{w(T-t)} - 1 \right),
			\end{aligned}
		\end{equation*}
		with \( G_1 := -2(\beta^2 + 2G_0) \), \( G_2 := 2\alpha + w - \frac{2(1-\gamma-\Phi)\beta \rho_1}{\Phi+\gamma} \), \( G_3 := - \frac{1-\gamma}{2 k(\Phi+\gamma)} \) and \( G_0 := \frac{(1-\gamma-\Phi)^2 k \beta^2 \rho_1^2}{2 (\Phi+\gamma) (1-\gamma)^2} - \frac{\Phi k^2 \beta^2}{2(1-\gamma)} + \frac{\beta^2 (k-1)}{2}\). The corresponding  strategies are
		\begin{equation}\label{piqc***CS}
			\begin{aligned}
				\pi^{\textnormal{CS}}(t) & = \frac{x(t)}{(\Phi+\gamma)\sigma^2} \left[(\sigma m(t)+a-r)  + \frac{1-\gamma-\Phi}{1-\gamma} k \beta \rho_1 \sigma \frac{g^{CS}_m(t, m(t))}{g^{CS}(t, m(t))} \right], \\
				q^{\textnormal{CS}}(t) & = \frac{\theta_1 \mu_1}{(\Phi+\gamma) \mu_2} x(t) , \\
				c^{\textnormal{CS}}(t) & =  x(t) \delta^\varphi (g^{CS}(t,m))^{k \frac{1-\varphi}{1-\gamma}},
			\end{aligned}
		\end{equation}
		and the distortion process \( \xi^{\textnormal{CS}} = \left(\xi_1^{\textnormal{CS}}, \xi_2^{\textnormal{CS}}, \xi_3^{\textnormal{CS}} \right)\) is
		\begin{equation*}\label{xiCS}
			\begin{aligned}
				\xi_1^{\textnormal{CS}}(t) & = \frac{\Phi}{(\Phi+\gamma)\sigma} (\sigma m(t)+a-r) + \frac{\Phi k \beta \rho_1}{(\Phi+\gamma)(1-\gamma)} \left(2 G(t) m(t) + L(t)\right), \\
				\xi_2^{\textnormal{CS}}(t) & = \frac{\Phi}{1-\gamma} k \beta \sqrt{1-\rho_1^2} \left(2 G(t) m(t) + L(t)\right), \\
				\xi_3^{\textnormal{CS}}(t) & = \frac{\Phi \theta_1 \mu_1 \sqrt{\lambda}}{(\Phi+\gamma) \sqrt{\mu_2}}.
			\end{aligned}
		\end{equation*}
	\end{thm}

	\begin{rem}
		Comparing the exact solution \eqref{piqc***} with the Campbell--Shiller approximation \eqref{piqc***CS} shows that the optimal reinsurance strategy  \( q^* \) is identical in both cases. 
		While the optimal investment--consumption strategies share the same functional forms, their specific values differ due to the function \( g \). Notably, the approximate optimal investment  ratio from the Campbell--Shiller method is independent of the reinsurance safety loading coefficient $\theta_1$, whereas the exact solution depends on $\theta_1$, and this discrepancy varies with  $\theta_1$.
	\end{rem}

	In the remainder of this section, we conduct a comparative analysis between the approximate solution and the exact solution. To demonstrate the flexibility of the model parameters, we adjust some of the baseline parameter settings in Table \ref{para}. Unless otherwise specified below, all other parameter settings remain consistent with those in Table \ref{para}.
	\begin{equation*}
		\gamma = 1.3, \quad \alpha = 7, \quad \Phi = 0.
	\end{equation*}

	\setlength{\abovecaptionskip}{2pt}
	\begin{figure}[h]
		\begin{minipage}[b]{0.49\linewidth}
			\centering
			\includegraphics[width=\textwidth]{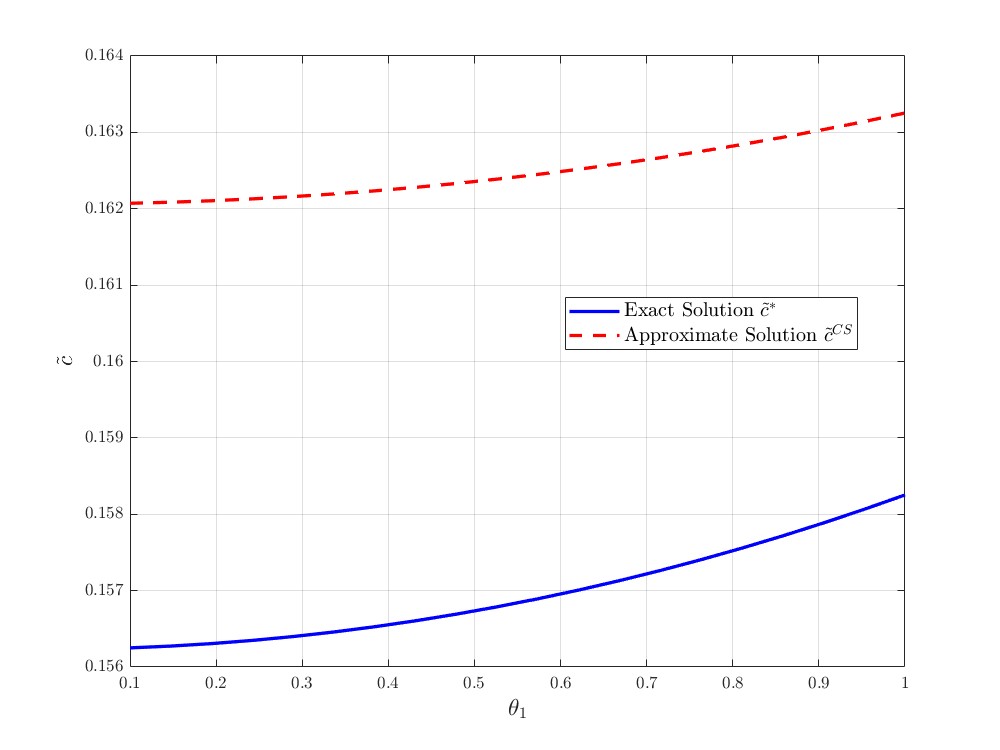}
			\centerline{\small(a) the reinsurance safety loading coefficient $\theta_1$}
		\end{minipage}
		\hfill
		\begin{minipage}[b]{0.49\linewidth}
			\centering
			\includegraphics[width=\textwidth]{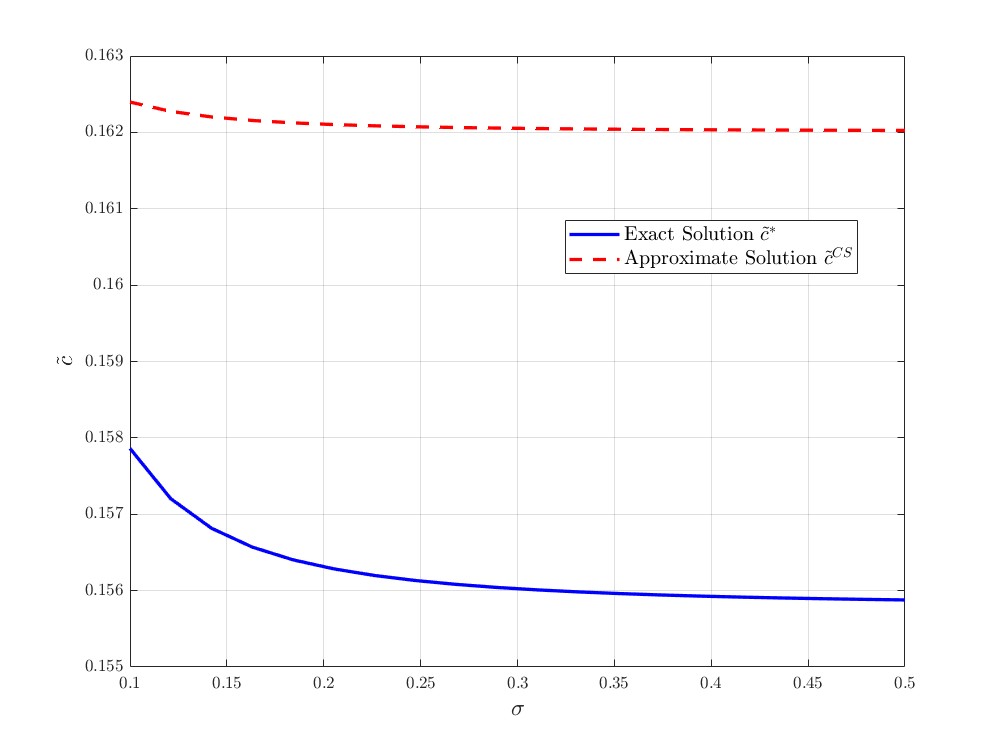}
			\centerline{\small(b) volatility of the risky asset $\sigma$}
		\end{minipage}
		
		\vspace{0.3cm}
		
		\caption{Comparison of $\tilde{c}^*$ and $\tilde{c}^{\text{CS}}$}
		\label{fig51}
	\end{figure}

	As shown in Figure \ref{fig51}, as the reinsurance safety loading coefficient $\theta_1$ increases or the volatility of the risky asset $\sigma$ decreases, the difference between the approximate solution and the exact solution gradually diminishes. Thus, when $\theta_1$ is sufficiently large or $\sigma$ is sufficiently small, the approximate solution provides a good approximation to the exact solution and can therefore serve as a substitute to simplify the analysis.

	\setlength{\abovecaptionskip}{2pt}
	\begin{figure}[h]
		\begin{minipage}[b]{0.49\linewidth}
			\centering
			\includegraphics[width=\textwidth]{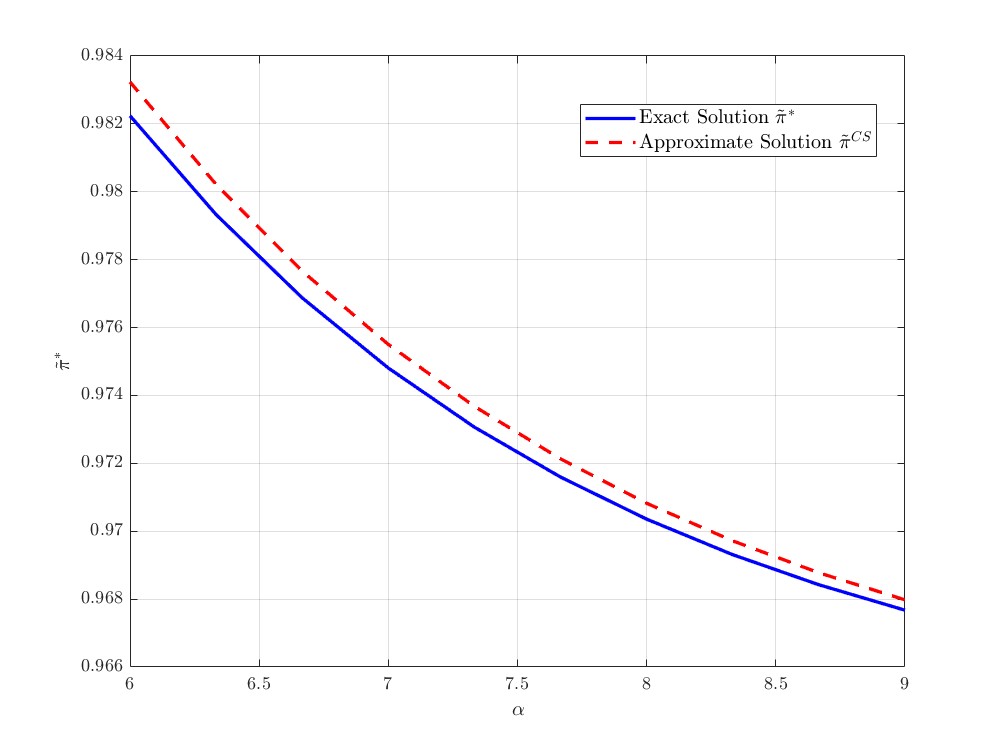}
			\centerline{\small(a) the mean reversion rate $\alpha$}
		\end{minipage}
		\hfill
		\begin{minipage}[b]{0.49\linewidth}
			\centering
			\includegraphics[width=\textwidth]{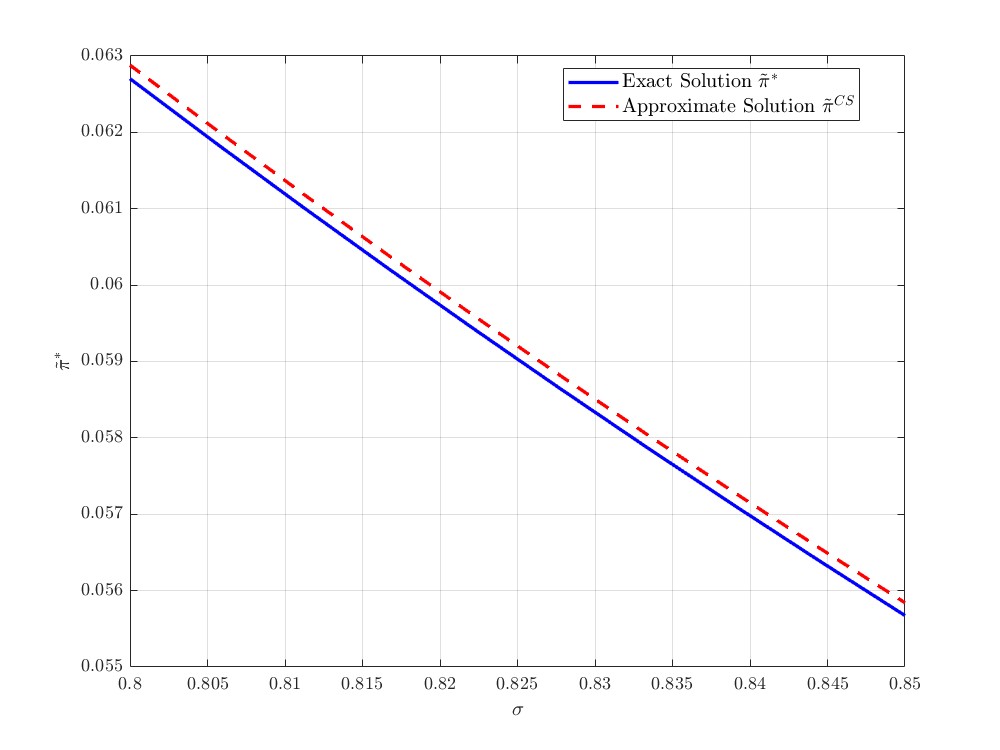}
			\centerline{\small(b) volatility of the risky asset $\sigma$}
		\end{minipage}
		
		\vspace{0.3cm}
		
		\caption{Comparison of $\tilde{\pi}^*$ and $\tilde{\pi}^{\text{CS}}$}
		\label{fig52}
	\end{figure}

	Besides, it can be observed from Figure \ref{fig52} and Table \ref{pisigma} that the difference between the approximate solution and the exact solution for the investment--wealth ratio is small under the following conditions: the mean-reversion speed $\alpha$ is small, and the volatility of the risky asset $\sigma$ is large. In these cases, the approximate solutions provide good approximations to the exact solutions, making it more reasonable to use them in practical applications.

	\renewcommand{\arraystretch}{2}
	\begin{table}[htbp]
		\scriptsize
		\centering
		\caption{Numerical results for $\pi$ with different $\sigma$}
		\label{pisigma}
		\begin{tabular}{c|c|c|c}
			\hline
			\textbf{$\sigma$} & \textbf{$\tilde{\pi}^{\text{CS}}$} & \textbf{$\tilde{\pi}^*$} & \textbf{Error} \\
			\hline
			0.8000 & 0.062880 & 0.062702 & 0.000178 \\
			0.8100 & 0.061368 & 0.061193 & 0.000175 \\
			0.8200 & 0.059911 & 0.059737 & 0.000173 \\
			0.8300 & 0.058505 & 0.058334 & 0.000171 \\
			0.8400 & 0.057150 & 0.056980 & 0.000169 \\
			0.8500 & 0.055841 & 0.055674 & 0.000167 \\
			\hline
		\end{tabular}
	\end{table}

	\section{Conclusion}\label{sec6}

	It is well known that the business objectives of an insurance company include not only external risk management (e.g., reinsurance) and investment appreciation, but also the ``internal consumption'' of its economic surplus, which is mainly manifested in returning capital to shareholders (such as dividend distributions and share repurchases) as well as performance-based employee compensation. 
	In this paper, taking internal consumption into account, we study the optimal robust reinsurance, investment, and surplus allocation problems for an insurance company under model uncertainty within the framework of Epstein--Zin recursive utility. Regarding model specification, we assume that the price of the risky asset follows a diffusion model with a stochastic drift driven by an O-U process, which better reflects the characteristics of actual financial markets. By applying the dynamic programming principle and solving the HJBI equation, we obtain explicit closed-form solutions for the optimal strategies and the corresponding value function. Since the introduction of the O-U process renders existing verification theorems inapplicable, we establish a new verification theorem by investigating the moment-generating properties of the O-U process and constructing an auxiliary process following a measure change. Unlike Han and Hung \cite{HH2017}, who provide only a single form of solution, for the case of non-unit EIS we provide both a Campbell--Shiller approximate solution (which is convenient for empirical applications) and a theoretically rigorous explicit closed-form solution, thereby significantly enhancing the model's applicability. Moreover, the robust decision-making framework constructed in this paper is highly general and compatible. When the intensity of model uncertainty approaches zero, our results degenerate to the classic non-robust case of Kraft et al. \cite{KSS2013}.

	\appendix
	\setcounter{equation}{0}                 
	\renewcommand{\theequation}{A.\arabic{equation}}

	\section*{Appendix}

	\subsection*{Proof of Proposition \ref{pro-exact}:}

	Based on the form given in \eqref{conject1}, we derive
	\begin{alignat*}{2}
		v_t    &= \frac{x^{1-\gamma}}{1-\gamma} k g^{k-1} g_t, \quad & v_x    &= x^{-\gamma} g^k, \\
		v_{xx} &= -\gamma x^{-\gamma-1} g^k, \quad & v_m    &= \frac{x^{1-\gamma}}{1-\gamma} k g^{k-1} g_m, \\
		v_{mm} &= \frac{x^{1-\gamma}}{1-\gamma} \left( k g^{k-1} g_{mm} + k(k - 1) g^{k-2} g_m^2 \right), \quad & v_{xm} &= x^{-\gamma} k g^{k-1} g_m.
	\end{alignat*}
	
	Substituting the derived derivatives into \eqref{piqc*} gives the candidate optimal strategies \eqref{piqc**}.
	Subsequently, by substituting both these derivatives and the candidate optimal solutions into the HJBI equation \eqref{HJBI2} and multiplying both sides by \( x^{\gamma-1} g^{-k} \), the HJBI equation \eqref{HJBI2} simplifies, after algebraic manipulation, to:
	\begin{equation}\label{HJBI3}
		\begin{aligned}
			& \left(1 - \frac{1}{\varphi}\right)^{-1} \delta \left(\delta^{\varphi-1} g^{k\frac{1-\varphi}{1-\gamma}} - 1\right) + \frac{(\sigma m+a-r)^2}{2(\Phi+\gamma)\sigma^2} + \frac{\lambda \theta_1^2 \mu_1^2}{2(\Phi+\gamma)\mu_2} + \frac{k}{1-\gamma} \frac{g_t}{g} + \left(r- \delta^\varphi g^{k\frac{1-\varphi}{1-\gamma}} \right) \\
			& + \frac{\beta^2 k}{2(1-\gamma)} \frac{g_{mm}}{g} + \left[\frac{1}{(\Phi+\gamma)\sigma^2} \left(\sigma m + a - r\right) k\beta \rho_1 \sigma \left(1-\frac{\Phi}{1-\gamma}\right) - \frac{\alpha m k}{1-\gamma} \right] \frac{g_m}{g} \\
			& + \left[\frac{1}{2(\Phi+\gamma)\sigma^2} k^2 \beta^2 \rho_1^2 \sigma^2 \left(1-\frac{\Phi}{1-\gamma} \right)^2 - \frac{\Phi k^2 \beta^2}{2(1-\gamma)^2} + \frac{\beta^2 k (k-1)}{2(1-\gamma)} \right]\left(\frac{g_m}{g}\right)^2 = 0.
		\end{aligned}
	\end{equation}
	Upon inserting the formula for \( k \) given in Assumption \ref{kvarphi} into the coefficient of \(\left( \frac{g_m}{g} \right)^2\), the coefficient reduces to zero. After this substitution, multiplying equation \eqref{HJBI3} by \( g(1-\gamma) k^{-1} \) yields the following HJB equation:
	\begin{equation}\label{HJBI4}
		\begin{aligned}
			& g_t + \frac{\frac{1}{\varphi}(1-\gamma)}{1-\frac{1}{\varphi}} \frac{\delta^\varphi}{k} g^{k \frac{1-\varphi}{1-\gamma} + 1} + \frac{1-\gamma}{k} \left[ r + \frac{(\sigma m + a - r)^2}{2(\Phi+\gamma)\sigma^2} + \frac{\lambda \theta_1^2 \mu_1^2}{2(\Phi+\gamma)\mu_2} - \frac{\delta}{1-\frac{1}{\varphi}} \right] g \\
			& + \left[ \frac{(1-\gamma-\Phi)}{(\Phi+\gamma)\sigma} (\sigma m + a - r) \beta \rho_1 - \alpha m \right] g_m + \frac{1}{2} \beta^2 g_{mm} = 0
		\end{aligned}
	\end{equation}
	Setting \( \varphi = 2 - \gamma - \Phi + \frac{(1 - \gamma - \Phi)^2}{\Phi + \gamma} \rho_1^2 \), which satisfies \( k \frac{1-\varphi}{1-\gamma} + 1 = 0 \), we find that equation \eqref{HJBI4} simplifies to equation \eqref{g1}. This completes the proof.
	\qed

	\subsection*{Proof of Theorm \ref{thm-exact}:}

	The proof is divided into six steps. Some of the proof ideas draw on the proofs in Appendix C of Kraft et al. \cite{KSS2013} and Theorem 4.2 of Pu and Zhang \cite{PZ2021}. In the specific derivation, technical details that are identical to those in the aforementioned two papers are omitted for brevity.

	\medskip
	\noindent \textbf{Step 1: Explicit solution of \eqref{g1}.} In this study, the PDE \eqref{h1} for \( h \) becomes
	\begin{equation}\label{h2}
		\begin{aligned}
			& h_t + \frac{1-\gamma}{k} \left[ r + \frac{(\sigma m + a - r)^2}{2(\Phi+\gamma)\sigma^2} + \frac{\lambda \theta_1^2 \mu_1^2}{2(\Phi+\gamma)\mu_2} - \frac{\delta}{1-\frac{1}{\varphi}} \right] h \\
			& + \left[ \frac{(1-\gamma-\Phi)}{(\Phi+\gamma)\sigma} (\sigma m + a - r) \beta \rho_1 - \alpha m \right] h_m + \frac{1}{2} \bar{\beta}^2 h_{mm} = 0.
		\end{aligned}
	\end{equation}
	
	We conjecture that the solution \( h \) has the following form:
	\begin{equation}\label{conjecth}
		h(t, m; s) = e^{A(t,s)-B(t,s)m-C(t,s)m^2}
	\end{equation}
	with \( h(s, m; s) = 1 \).
	Substituting \eqref{conjecth} into \eqref{h2}, and classify the terms involving \(m^2, m\), and those without \(m\). Then setting the coefficients of each term to zero yields the following three ordinary differential equations (ODEs):
	\begin{equation}\label{ABCode}
		\begin{aligned}
			& C'(t, s) + 2 \left( \frac{1-\gamma-\Phi}{\Phi+\gamma}\beta \rho_1 - \alpha \right) C(t, s) - 2\beta^2 C^2(t, s) - \frac{1-\gamma}{2 k (\Phi+\gamma)} = 0, \quad C(s, s) = 0;  \\
			& B'(t, s) + \left[ \frac{1-\gamma-\Phi}{\Phi+\gamma}\beta \rho_1 - \alpha - 2\beta^2 C^2(t, s) \right] B(t, s) + \frac{1-\gamma-\Phi}{(\Phi+\gamma)\sigma} \beta \rho_1 (a-r) C(t, s) \\
			& - \frac{(1-\gamma)(a-r)}{k (\Phi+\gamma)\sigma} = 0, \quad B(s, s) = 0;  \\
			& A'(t, s) + \frac{1}{2} \beta^2 B^2(t,s) - \beta^2 C(t,s) - \frac{1-\gamma-\Phi}{(\Phi+\gamma)\sigma} (a-r) \beta \rho_1 B(t,s) + \frac{1-\gamma}{k} \Bigg[ r + \frac{(a-r)^2}{2(\Phi+\gamma)\sigma^2} \\
			&  + \frac{\lambda \theta_1^2 \mu_1^2}{2(\Phi+\gamma)\mu_2} - \frac{\delta}{1-\frac{1}{\varphi}} \Bigg] = 0, \quad A(s, s) = 0.
		\end{aligned}
	\end{equation}
	It is easy to get that \eqref{ABCsolution} solves above ODEs \eqref{ABCode}.
	The candidate optimal strategy \eqref{piqc***} follows from the results \eqref{piqc**}.

	\medskip
	\noindent \textbf{Step 2: Verification that \((\pi^*, q^*, c^*; \xi^*)\) is admissible.} Prior to verifying that the optimal strategy is admissible, we first provide the bounds for \(A, B\) and \( C \).
	Based on the expressions for \(A, B\) and \( C \) in \eqref{ABCsolution} and the standard properties of ODEs, their respective bounds are derived as:
	\begin{equation}\label{ABCbound}
		\begin{aligned}
			& 0 \leq C(t,s) \leq b_0(s-t), \ C(t,s) \leq \frac{2b_0}{2\kappa+\Delta}, \\
			& 0 \leq |B(t,s)| \leq b_1 (s-t), \\
			& A(t,s) \geq A_1(T-t)(s-t) + A_2(s-t),
		\end{aligned}
	\end{equation}
	where
	\begin{equation*}
		\begin{aligned}
			& A_1 = - \frac{(1-\gamma-\Phi)(1-\gamma)(a-r)^2 \beta \rho_1 }{k(\Phi+\gamma)^2\sigma^2} \left( \frac{4(1-\gamma-\Phi) b_0 k \beta \rho_1}{(1-\gamma) (2\kappa + \Delta)} - 1\right) ,\\
			& A_2 = \frac{1-\gamma}{k} \left( r + \frac{(a-r)^2}{2(\Phi+\gamma)\sigma^2} + \frac{\lambda \theta_1^2 \mu_1^2}{2(\Phi+\gamma) \mu_2} \right) - \frac{2b_0 \beta^2}{2\kappa+\Delta}, \\
			& b_1 = \frac{(1-\gamma)(a-r)}{k(\Phi+\gamma)\sigma} \left( \frac{4(1-\gamma-\Phi) b_0 k \beta \rho_1}{(1-\gamma) (2\kappa + \Delta)} - 1 \right).
		\end{aligned}
	\end{equation*}
	Under the parameter settings previously defined, it follows directly that \(A_1\) and \(A_2\) are negative.
	For \(c^*\), let
	\begin{equation*}
		\tilde{c}^* := \frac{c^*}{x} = \frac{\delta^\varphi}{g} = \frac{1}{\int_t^T h(t, m; s) \, \textnormal{d}s + \delta^{-\varphi} h(t, m; T)}  \leq \frac{1}{\int_t^T h(t, m; s) \, \textnormal{d}s},
	\end{equation*}
	and from \eqref{ABCbound} we have
	\begin{align}
		\int_t^T h(t, m; s) \, \textnormal{d}s
		& = \int_t^T e^{A(t,s) - B(t,s)m - C(t,s)m^2}  \, \textnormal{d}s \nonumber \\
		& \geq \int_t^T e^{A(t,s) - |B(t,s)||m| - C(t,s)m^2}  \, \textnormal{d}s \nonumber \\
		& \geq \int_t^T e^{A_1(T-t)(s-t) + A_2(s-t) - b_1(s-t)|m| - b_0(s-t)m^2}  \, \textnormal{d}s \nonumber \\
		& = \int_t^T e^{-[b_0 m^2 + b_1|m| - A_2 - A_1(T-t)](s-t)} \, \textnormal{d}s \nonumber \\
		& = \frac{1 - e^{-[b_0 m^2 + b_1|m| - A_2 - A_1(T-t)](T-t)}}{b_0 m^2 + b_1|m| - A_2 - A_1(T-t)} \nonumber \\
		& = \frac{1}{\tilde{h}(b_0 m^2 + b_1|m| - A_2 - A_1(T-t))} \nonumber,
	\end{align}
	where \( \tilde{h}(x) = \frac{x}{1 - e^{x(T-t)}}\). By using (C.7) in Lemma C.5 of Kraft et al. \cite{KSS2013}, we obtain
	\begin{equation*}
		\frac{1}{\int_t^T h(t, m; s) \, \textnormal{d} s} \leq \frac{1}{T-t} + b_0 m^2 + b_1|m| - A_2 - A_1(T-t),
	\end{equation*}
	which indicates \( c^* \in \mathcal{A}_t \).

	For \(\pi^*\), we first investigate the upper bound of \( |g_m/g| \), where \(g\) is defined in \eqref{g*}.
	Decompose the interval \([t,T]\) into \([t,t+\varepsilon]\) and \([t+\varepsilon, T]\), and let
	\begin{equation*}
		\begin{aligned}
			& I_1(t,m) = \int_t^{t+\varepsilon}  h(t, m; s) \textnormal{d}s, \ J_1(t,m) = \int_t^{t+\varepsilon}  h(t, m; s) (-B(t,s)-2C(t,s)m) \textnormal{d}s, \\
			& I_2(t,m) = \int_{t+\varepsilon}^{T}  h(t, m; s) \textnormal{d}s, \ J_2(t,m) = \int_{t+\varepsilon}^{T}  h(t, m; s) (-B(t,s)-2C(t,s)m) \textnormal{d}s.
		\end{aligned}
	\end{equation*}
	Then, \( |g_m/g| \) can be simplified to
	\begin{equation*}
		\begin{aligned}
			\frac{g_m}{g} & = \frac{\delta^\varphi \int_t^T h(t, m; s) (-B(t,s)-2C(t,s)m) \textnormal{d}s +  h(t, m; T) (-B(t,T)-2C(t,T)m) }{\delta^\varphi \int_t^T  h(t, m; s) \textnormal{d}s +  h(t, m; T)} \\
			& = \frac{\delta^\varphi (J_1(t,m) + J_2(t,m)) +  h(t, m; T) (-B(t,T)-2C(t,T)m)}{\delta^\varphi (I_1(t,m) + I_2(t,m)) +  h(t, m; T)}.
		\end{aligned}
	\end{equation*}
	By analyzing the asymptotic behavior on each subinterval as \( |m| \rightarrow \infty \), we arrive at
	\begin{equation}\label{I1J1}
		I_1(t,m) = \frac{1}{K_t m^2} + o\left(\frac{1}{m^2}\right), \
		J_1(t,m) = - \frac{2}{K_t |m|^3} + o\left(\frac{1}{|m|^3}\right),
	\end{equation}
	where \(K_t\) is the value of the first derivative of \( C(t,s) \)  with respect to \( s \) evaluated at \( s=t \), which is computed to be \( b_0 \) and is independent of the choice of the initial time \( t \).
	As for \( I_2(t,m) \), \( J_2(t,m) \), and the terminal terms \( h(t, m; T) \) and \( h(t, m; T) \bigl( -B(t, T) - 2C(t, T)m \bigr) \), when \( |m| \) is sufficiently large, their asymptotic behavior is dominated by \( e^{- (c_\varepsilon / 2) m^2} \), where \( c_\varepsilon = \min_{s \in [t+\varepsilon, T]} C(t,s) > 0 \). Consequently, these terms are negligible compared to \( I_1(t,m) \) and \( J_1(t,m) \). 
	Therefore,
	\begin{equation}\label{gm/g}
		\begin{aligned}
			\frac{g_m}{g} & = \frac{\delta^\varphi (J_1(t,m) + J_2(t,m)) +  h(t, m; T) (-B(t,T)-2C(t,T)m)}{\delta^\varphi (I_1(t,m) + I_2(t,m)) +  h(t, m; T)} \\
			& = \frac{J_1(t,m)}{I_1(t,m)}(1+o(1)) = \frac{A_g}{|m|} + o\left(\frac{1}{|m|}\right),
		\end{aligned}
	\end{equation}
	where \(A_g\) is a constant, which is independent of \(t\).
	Hence, there exists a constant \( K_g > 0 \) such that for all sufficiently large \( |m| \), the following inequality holds:
	\begin{equation*}
		|\frac{g_m}{g}| \leq \frac{K_g}{1 + |m|}, \quad \forall t \in [0, T], \, m \in \mathbb{R}.
	\end{equation*}
	Consequently, the following inequality holds for \( \tilde{\pi}^* \):
	\begin{equation*}
		|\tilde{\pi}^*| = |\frac{\pi^*}{x}| \leq \frac{1}{(\Phi+\gamma)\sigma} m +  \frac{1}{(\Phi+\gamma)\sigma^2} \left[ a - r  + \left(1-\frac{\Phi}{1-\gamma}\right) k \beta \rho_1 \sigma K_g \right],
	\end{equation*}
	which shows that \( \pi^* \in \mathcal{A}_t \).
	
	The expression for \( q^* \) readily confirms that \( q^* \in \mathcal{A}_t \).
	
	Finally, for the distortion process \(\xi^*\) \eqref{xi**}, the verification method for \(\xi_1^*\) follows a similar approach to that of \(\pi^*\), while the verification processes for \(\xi_2^*\) and \(\xi_3^*\) are relatively straightforward. Therefore, the conclusion that \( \xi^* \in \mathcal{A}_t \) follows directly from the preceding analysis and will not be reiterated here.

	\medskip
	\noindent \textbf{Step 3: Verification of the condition that the local martingale
		$$
		\int_t^\cdot v_x \sigma \pi(s) \textnormal{d}W_{1,s}^{\mathbb{Q}^\xi} + \int_t^\cdot v_x \sqrt{\lambda \mu_2} q(s) \textnormal{d}W_{3,s}^{\mathbb{Q}^\xi} + \int_t^\cdot v_m \beta\left( \rho_1 + \sqrt{1-\rho_1^2} \right) \textnormal{d}W_{2,s}^{\mathbb{Q}^\xi}
		$$
		is a true martingale for any \( (\pi, q, c; \xi) \in \mathcal{A}_t \).}
	For the above condition to hold, it suffices to satisfy the following integrability condition:
	\begin{equation}\label{keji}
		\mathbb{E}^{\mathbb{Q}^\xi} \left\{ \int_0^T \left[ (v_x \sigma \pi(t))^2 + (v_x \sqrt{\lambda \mu_2} q(t))^2 + \left( v_m \beta\left( \rho_1 + \sqrt{1-\rho_1^2} \right) \right)^2 \right] \textnormal{d}t \right\} < \infty.
	\end{equation}
	We analyze each term in \eqref{keji}. Substituting the strategy \eqref{piqc***} and the expressions for \( v_x \) and \( v_m \), the three integrands can be written as
	\begin{equation*}
		\begin{aligned}
			& v_x \sigma \pi(t) = \frac{X_t^{1-\gamma} \left( g(t, m) \right)^k}{(\Phi + \gamma)\sigma} \left[ (\sigma m(t) + a - r) + \left(1 - \frac{\Phi}{1-\gamma}\right) k \beta \rho_1 \sigma \frac{g_m(t, m)}{g(t, m)} \right], \\
			& v_x \sqrt{\lambda \mu_2} q(t) = \frac{\sqrt{\lambda \mu_2} \theta_1 \mu_1}{(\Phi+\gamma) \mu_2} X_t^{1-\gamma} \left( g(t, m) \right)^k, \\
			& v_m \beta\left( \rho_1 + \sqrt{1-\rho_1^2} \right) = \beta\left( \rho_1 + \sqrt{1-\rho_1^2} \right) \frac{k}{1-\gamma} X_t^{1-\gamma} \left( g(t, m) \right)^k \frac{g_m(t, m)}{g(t, m)}.
		\end{aligned}
	\end{equation*}
	Based on the above expression, it is necessary to examine the boundedness of \( g(t,m) \) and its \( k \)-th power \( (g(t,m))^k \).
	Similar to the analysis of \( g_m/g \) in Step 2, we denote \(g(t,m)\) \eqref{g*} as \( g(t,m) = \delta^\varphi (I_1(t,m) + I_2(t,m)) + h(t, m; T) \). For \( I_2(t,m) \), let \( s = t+u \), and we obtain
	\begin{equation*}
		\begin{aligned}
			I_2(t,m) = \delta^\varphi \int_\varepsilon^{T-t} e^{A(t,t+u) - B(t,t+u)m - C(t,t+u)m^2} \, \textnormal{d}u \leq \delta^\varphi \int_\varepsilon^{T-t} e^{A_{max} + b_1 u |m| - c_1 u m^2} \, \textnormal{d}u,
		\end{aligned}
	\end{equation*}
	where \(A_{max} = \max_{u \in [0, T-t]} A(t,t+u) \), and \( c_1 = \min_{u \in [\varepsilon, T-t]} C(t,t+u)/u \). When \( |m| \) is sufficiently large, we have \( b_1 |m| - c_1 m^2 < 0 \). Consequently, the above inequality can be further written as
	\(
	I_2(t,m) \leq \delta^\varphi e^{A_{max}} (T-t).
	\)
	For \( I_1 \), when \( |m| \) is sufficiently large, it follows from \eqref{I1J1} that it tends to zero. For the terminal term \( h(t, m; T) \), it is straightforward to see that it also tends to zero. Therefore, when \( |m| \) is sufficiently large, \( g(t,m) \) satisfies the following inequality:
	\begin{equation*}
		0 \leq g(t, m) \leq \delta^\varphi e^{A_{max}} (T - t) \triangleq M_g (T-t), \quad t \in [0, T], \, m \in \mathbb{R}.
	\end{equation*}
	In addition, the above inequality remains valid when \( |m| \)  is bounded.
	Combined with the bound \eqref{gm/g} for \(|g_m/g|\) when \(|m|\) is sufficiently large, as derived in Step 2, we have
	\begin{equation*}
		v_x \sigma \pi(t) \leq \frac{M_g^k}{(\Phi + \gamma)\sigma} (T - t)^k X_t^{1-\gamma} \left[ \sigma m(t) + (a - r) + \left(1 - \frac{\Phi}{1-\gamma}\right) k \beta \rho_1 \sigma K_g \right].
	\end{equation*}
	
	Let \( K_{g_1} = \frac{M_g^k}{(\Phi + \gamma)\sigma}\), \( K_{g_2} = a - r + \left(1 - \frac{\Phi}{1-\gamma}\right) k \beta \rho_1 \sigma K_g\), then
	\begin{equation*}
		(v_x \sigma \pi(t))^2 \leq K_{g_1}^2 (T-t)^{2k} X_t^{2(1-\gamma)}  (\sigma m(t) + K_{g_2})^2 \leq 2 K_{g_1}^2 (T-t)^{2k} X_t^{2(1-\gamma)}  (\sigma^2 m^2(t) + K_{g_2}^2).
	\end{equation*}
	Furthermore, as demonstrated in the proof of Lemma C.3 in Kraft et al. \cite{KSS2013}, the admissibility condition (M) implies that
	\begin{equation*}
		\mathbb{E}^{\mathbb{Q}^\xi} [X_t^{\overline{k}(1-\gamma)}] \leq K(T - t)^{\overline{k}(1-\gamma)}.
	\end{equation*}
	if \(\overline{k} > 2 \), and \(\overline{k} - 2\) is sufficiently small.
	By H$\ddot{o}$lder's inequality,
	\begin{equation*}
		\begin{aligned}
			\mathbb{E}^{\mathbb{Q}^\xi}[|X_t^{2(1-\gamma)} m^2(t)|]
			& \leq [\mathbb{E}^{\mathbb{Q}^\xi}|X_t^{2p(1-\gamma)}|]^{1/p} \cdot [\mathbb{E}^{\mathbb{Q}^\xi}|(m(t))^{2q}|]^{1/q} \leq K_m \cdot [\mathbb{E}^{\mathbb{Q}^\xi}|X_t^{2p(1-\gamma)}|]^{1/p} \\
			& \leq  K_m \cdot K (T-t)^{2(1-\gamma)},
		\end{aligned}
	\end{equation*}
	where \( K_m \) denotes the upper bound of the \( 2q \)-th moment of the new process \eqref{ROU} obtained from the O-U process after a measure transformation (based on the properties of the O-U process under measure transformation, it is not difficult to verify that the moments of all orders of the transformed new process remain finite).
	By (H1), \( 2(k+1-\gamma) >-1\). Consequently,
	\begin{equation*}
		\begin{aligned}
			\mathbb{E}^{\mathbb{Q}^\xi}\left[\int_0^T (v_x \sigma \pi(t))^2 \, \textnormal{d}t \right] & = \int_0^T \mathbb{E}^{\mathbb{Q}^\xi}[(v_x \sigma \pi(t))^2] \, \textnormal{d}t \\
			& \leq 2 K_{g_1}^2 \sigma^2 \int_0^T (T-t)^{2k} \mathbb{E}^{\mathbb{Q}^\xi} \left[ X_t^{2(1-\gamma)} m^2(t) \right] \, \textnormal{d}t \\
			& \quad + 2 K_{g_1}^2 K_{g_2}^2 \int_0^T (T-t)^{2k} \mathbb{E}^{\mathbb{Q}^\xi} \left[ X_t^{2(1-\gamma)} \right] \, \textnormal{d}t \\
			& \leq 2 K_{g_1}^2 K (\sigma^2 K_m + K_{g_2}^2) \int_0^T (T-t)^{2(k+1-\gamma)} \, \textnormal{d}t < +\infty.
		\end{aligned}
	\end{equation*}
	Similarly, we can also obtain
	\begin{equation*}
		\begin{aligned}
			& \mathbb{E}^{\mathbb{Q}^\xi}\left[\int_0^T (v_x \sqrt{\lambda \mu_2} q(t))^2 \, \textnormal{d}t \right] \leq \frac{\lambda \theta_1^2 \mu_1^2}{(\Phi+\gamma)^2 \mu_2} M_g^2 K \int_0^T (T-t)^{2(k+1-\gamma)}  \, \textnormal{d}t < +\infty, \\
			& \mathbb{E}^{\mathbb{Q}^\xi}\left[\int_0^T \left( v_m \beta  \left( \rho_1 + \sqrt{1-\rho_1^2} \right) \right)^2 \, \textnormal{d}t \right] \leq \frac{\beta^2 \left( \rho_1 + \sqrt{1-\rho_1^2} \right) k^2 K_g^2}{(1-\gamma)^2} M_g^2 K \cdot \\
			& \qquad\qquad\qquad\qquad\qquad\qquad\qquad\qquad\qquad \int_0^T (T-t)^{2(k+1-\gamma)}  \, \textnormal{d}t < +\infty,
		\end{aligned}
	\end{equation*}
	which further proves that the integrability condition \eqref{keji} holds.

	\medskip
	\noindent \textbf{Step 4: Verification of the condition (M).}
	Within the admissible control set \(\mathcal{A}_t\), we rewrite the wealth process \(X_t\) \eqref{RX} in terms of \(\tilde{\pi}\), \(\tilde{q}\), and \(\tilde{c}\), from which \(X_t\) can be computed explicitly as
	\begin{equation*}
		\begin{aligned}
			X_t = x \cdot \exp \Bigg\{ \int_0^t & \Bigg[ r + (\sigma m(s) + a - r)\tilde{\pi}(s) + \lambda \theta_1\mu_1 \tilde{q}(s) - \tilde{c}(s) - \sigma \xi_1^*(s) \tilde{\pi}(s) -  \sqrt{\lambda\mu_2} \xi_3^*(s) \tilde{q}(s) \\
			& - \frac{1}{2} \left( \sigma^2 \tilde{\pi}^2(s) + \lambda \mu_2 \tilde{q}^2(s) \right) \Bigg] \textnormal{d}s + \int_0^t \sigma \tilde{\pi}(s) \textnormal{d}W_{1,s}^{\mathbb{Q}^\xi}  + \int_0^t \sqrt{\lambda \mu_2}\tilde{q}(s) \textnormal{d}W_{3,s}^{\mathbb{Q}^\xi} \Bigg\}.
		\end{aligned}
	\end{equation*}
	It follows that
	\begin{equation*}
		X_t^{\overline{k}(1-\gamma)} \leq x^{\overline{k}(1-\gamma)} \exp \left\{ \overline{k}(\gamma-1) \int_0^t Y_s \, \textnormal{d}s \right\} Z_t,
	\end{equation*}
	where the processes \(Y = \{Y_t\} \) and \(Z = \{Z_t\} \) are given by
	\begin{equation*}
		\begin{aligned}
			& Y_t := - \Big[ r + (\sigma m(s) + a - r)\tilde{\pi}(t) + \lambda \theta_1\mu_1 \tilde{q}(t) - \tilde{c}(t) - \sigma \xi_1^*(t) \tilde{\pi}(t) - \sqrt{\lambda\mu_2} \xi_3^*(t) \tilde{q}(t) \\
			& \qquad\qquad - \frac{1}{2} ( \overline{k}(\gamma-1) + 1 ) \left( \sigma^2 \tilde{\pi}^2(t) + \lambda \mu_2 \tilde{q}^2(t) \right) \Big], \\
			& Z_t := \exp \Bigg\{ \overline{k}(1-\gamma) \int_0^t \sigma \tilde{\pi}(s) \textnormal{d}W_{1,s}^{\mathbb{Q}^\xi} + \overline{k}(1-\gamma) \int_0^t \sqrt{\lambda \mu_2}\tilde{q}(s) \textnormal{d}W_{3,s}^{\mathbb{Q}^\xi} - \frac{1}{2} \overline{k}^2(1-\gamma)^2 \\
			& \qquad\qquad\quad \int_0^t \left( \sigma^2 \tilde{\pi}^2(s) + \lambda \mu_2 \tilde{q}^2(s) \right) \textnormal{d}s \Bigg\}.
		\end{aligned}
	\end{equation*}
	By Novikov's condition and \textbf{Step 5}, \( Z = \{Z_t\} \) is a martingale. Hence, under the equivalent measure \( \tilde{\mathbb{Q}}^\xi \) given by
	\(
	\frac{d\tilde{\mathbb{Q}}^\xi}{d\mathbb{Q}^\xi} = Z_t
	\)
	on \( \mathcal{F}_t \), the process \( W_{1,t}^{\tilde{\mathbb{Q}}^\xi} \) and \( W_{3,t}^{\tilde{\mathbb{Q}}^\xi} \) defined by
	\begin{equation*}
		\begin{aligned}
			W_{1,t}^{\tilde{\mathbb{Q}}^\xi} := W_{1,t}^{\mathbb{Q}^\xi} - \overline{k}(1 - \gamma) \int_0^t \sigma \pi(s) \, \textnormal{d}s, \ \
			W_{3,t}^{\tilde{\mathbb{Q}}^\xi} := W_{3,t}^{\mathbb{Q}^\xi} - \overline{k}(1 - \gamma) \int_0^t \sqrt{\lambda \mu_2} q(s) \, \textnormal{d}s
		\end{aligned}
	\end{equation*}
	is a standard Wiener process. With \( \mathbb{E}^{\tilde{\mathbb{Q}}^\xi} \) denoting the expectation operator associated to \( \tilde{\mathbb{Q}}^\xi \), we thus have
	\begin{equation*}
		\mathbb{E}^{\mathbb{Q}^\xi} [X_t^{\overline{k}(1-\gamma)}] \leq x^{\overline{k}(1-\gamma)} \mathbb{E}^{\tilde{\mathbb{Q}}^\xi} \left[ \exp \left\{ \overline{k}(\gamma - 1) \int_0^t Y_s \, \textnormal{d}s \right\} \right].
	\end{equation*}
	
	Furthermore, under the measure \(\tilde{\mathbb{Q}}^\xi\), the original O-U process \(m(t)\) \eqref{303} evolves as
	\begin{equation}\label{m(t)trans2}
		\textnormal{d} m(t)
		= \left(-\bar{\alpha}_0 m(t) - \bar{\alpha}_1 - \bar{\alpha}_2 \frac{g_m(t,m)}{g(t,m)} \right)\textnormal{d}t + \beta \left[ \rho_1 \textnormal{d} W_{1,t}^{\tilde{\mathbb{Q}}^\xi} + \sqrt{1-\rho_1^2} W_{2,t}^{\tilde{\mathbb{Q}}^\xi} \right],
	\end{equation}
	where \( \bar{\alpha}_0 = \alpha + \frac{(\Phi + \overline{k}(\gamma-1)) \beta \rho_1}{\Phi+\gamma} > 0 \), \( \bar{\alpha}_1 = \frac{( \Phi + \overline{k}(\gamma-1) ) (a-r) \beta \rho_1 }{(\Phi+\gamma)\sigma} \), \( \bar{\alpha}_2 = \frac{(\Phi + \overline{k}(\gamma-1)) (1-\gamma-\Phi) k \beta^2 \rho_1^2 }{(\Phi+\gamma)(1-\gamma)} + \frac{\Phi k \beta^2}{1-\gamma} \), and \( W_{2,t}^{\tilde{\mathbb{Q}}^\xi} = W_{2,t}^{\mathbb{Q}^\xi} \).
	According to the properties of the O-U process and the measure transformation, the moments of all orders of \( m(t) \) \eqref{m(t)trans2} under the measure \(\tilde{\mathbb{Q}}^\xi\) are finite.
	
	By virtue of the conditions of the admissible control set \(\mathcal{A}_t\), we scale \( Y_t \) as follows:
	\begin{equation*}
		Y_t \leq D_7 m^2(t) + D_8 |m(t)| + D_9 + \frac{1}{T-t},
	\end{equation*}
	where
	\begin{equation*}
		\begin{aligned}
			& D_7 := D_1 + \sigma D_1D_4 + b_0 + \frac{1}{2} (\overline{k}(\gamma-1) + 1) \sigma^2 D_1^2, \\
			& D_8 := D_1(a-r) + D_2 + \sigma (D_1D_5 + D_2D_4) + b_1 + (\overline{k}(\gamma-1) + 1) \sigma^2 D_1^2 , \\
			& D_9 := D_2(a-r) + \sigma D_2 D_5 + \sqrt{\lambda\mu_2} D_3 D_6 + \frac{1}{2} (\overline{k}(\gamma-1) + 1) (\sigma^2 D_2^2 + \lambda \mu_2 D_3^2) - A_1T - A_2 .
		\end{aligned}
	\end{equation*}
	Consequently, we have
	\begin{equation}\label{EYt}
		\begin{aligned}
			\mathbb{E}^{\mathbb{Q}^\xi} [X_t^{\overline{k}(1-\gamma)}] & \leq x^{\overline{k}(1-\gamma)} \mathbb{E}^{\tilde{\mathbb{Q}}^\xi} \left[ \exp \left\{ \overline{k}(\gamma - 1) \int_0^t \left( D_7 m^2(s) + D_8 |m(s)| + D_9 + \frac{1}{T-t} \right) \, \textnormal{d}s \right\} \right] \\
			& \leq x^{\overline{k}(1-\gamma)} \exp \left\{ \overline{k}(\gamma - 1)  D_9 T \right\} T^{\overline{k}(\gamma - 1)} (T-t)^{-\overline{k}(\gamma - 1)} \mathbb{E}^{\tilde{\mathbb{Q}}^\xi} \Big[ \exp \Big\{ \overline{k}(\gamma - 1) D_7 \\
			& \quad \  \int_0^t m^2(s) \, \textnormal{d}s \Big\} \cdot \exp \Big\{ \overline{k}(\gamma - 1) D_8 \int_0^t |m(s)| \, \textnormal{d}s \Big\} \Big] \\
			& \leq x^{\overline{k}(1-\gamma)} \exp \left\{ \overline{k}(\gamma - 1)  D_9 T \right\} T^{\overline{k}(\gamma - 1)} (T-t)^{-\overline{k}(\gamma - 1)} Y_1^{1/2}(t) \cdot Y_2^{1/2}(t),
		\end{aligned}
	\end{equation}
	where
	\begin{equation*}
		\begin{aligned}
			& Y_1(t) := \mathbb{E}^{\tilde{\mathbb{Q}}^\xi} \left[ \exp \left\{ 2 \overline{k}(\gamma - 1) D_7 \int_0^t m^2(s) \, \textnormal{d}s \right\} \right], \\
			& Y_2(t) := \mathbb{E}^{\tilde{\mathbb{Q}}^\xi} \left[ \exp \left\{ 2 \overline{k}(\gamma - 1) D_8 \int_0^t |m(s)| \, \textnormal{d}s \right\} \right].
		\end{aligned}
	\end{equation*}

	For \(Y_2(t)\), we construct an auxiliary O-U process
	\begin{equation}\label{auxiOU}
		\textnormal{d}\bar{m}(t) = -\bar{\alpha}_0 \bar{m}(t) \textnormal{d}t + \beta \left[ \rho_1 \textnormal{d} W_{1,t}^{\tilde{\mathbb{Q}}^\xi} + \sqrt{1-\rho_1^2} \textnormal{d} W_{2,t}^{\tilde{\mathbb{Q}}^\xi} \right], \quad \bar{m}(0) = m.
	\end{equation}
	Define \(\bar{\Delta}(t) = m(t) - \bar{m}(t)\). Then
	\begin{equation*}\label{auxiDelta1}
		\textnormal{d}\bar{\Delta}(t) = -\bar{\alpha}_0 \bar{\Delta}(t) \textnormal{d}t + (-\bar{\alpha}_1 - \bar{\alpha}_2 \frac{g_m}{g}) \textnormal{d}t, \quad \bar{\Delta}(0) = 0.
	\end{equation*}
	Furthermore,
	\begin{equation*}
		\bar{\Delta}(t) = \int_0^t e^{-\bar{\alpha}_0 (t-s)} \left( -\bar{\alpha}_1 - \bar{\alpha}_2 \frac{g_m}{g} \right) \textnormal{d}s.
	\end{equation*}
	Since \(|-\bar{\alpha}_1 - \bar{\alpha}_2 \frac{g_m}{g}| \leq D_{10}\), where \(D_{10} = |\bar{\alpha}_1| + |\bar{\alpha}_2| K_g \), it follows that
	\begin{equation*}\label{auxiDelta2}
		|\bar{\Delta}(t)| \leq D_{10} \int_0^t e^{-\bar{\alpha}_0 (t-s)} \textnormal{d}s = \frac{D_{10}}{\bar{\alpha}_0} (1-e^{-\bar{\alpha}_0 t}) \leq \frac{D_{10}}{\bar{\alpha}_0}.
	\end{equation*}
	Thus,
	\begin{equation*}
		|m(t)| = |\bar{m}(t) + \bar{\Delta}(t)| \leq |\bar{m}(t)| + |\bar{\Delta}(t)| \leq |\bar{m}(t)| + \frac{D_{10}}{\bar{\alpha}_0}.
	\end{equation*}

	Based on the solution form and properties of the O-U process, we have
	\[
	\bar{m}(t) = m e^{-\bar{\alpha}_0 t} + \beta \int_{0}^{t} e^{-\bar{\alpha}_0 (t-s)} \, \textnormal{d}W_s^{\tilde{\mathbb{Q}}^\xi} := \bar{\mu}(t) + \bar{Z}(t), \ \
	\bar{m}(t) \sim N(\bar{\mu}(t), \bar{\sigma}^2(t)),
	\]
	where
	\(
	\textnormal{d}W_s^{\tilde{\mathbb{Q}}^\xi} = \rho_1 \, \textnormal{d}W_{1,s}^{\tilde{\mathbb{Q}}^\xi} + \sqrt{1 - \rho_1^2} \, \textnormal{d}W_{2,s}^{\tilde{\mathbb{Q}}^\xi},
	\)
	\( \bar{\sigma}^2(t) = \frac{\beta^2}{2\bar{\alpha}_0} (1 - e^{-2\bar{\alpha}_0 t}) \),
	and \(\bar{Z}(t)\) is a Gaussian process with zero mean and variance 
	\( \operatorname{Var}(\bar{Z}(t)) = \bar{\sigma}^2(t) \leq \frac{\beta^2}{2\bar{\alpha}_0}. \)
	Then
	\begin{equation*}
		\int_{0}^{t} |\bar{m}(s)| \, \textnormal{d}s \leq \int_{0}^{t} |\bar{\mu}(s)| \, \textnormal{d}s + \int_{0}^{t} |\bar{Z}(s)| \, \textnormal{d}s \leq m t + \sup_{s \in [0,t]} |\bar{Z}(s)| t.
	\end{equation*}
	Applying the Borell-Sudakov-Tsirelson concentration inequality (Theorem 2.2.7 in Gin$\acute{e}$ and Nickl \cite{GN2016}), we obtain
	\begin{equation}\label{BST1}
		\tilde{\mathbb{Q}}^\xi \left\{ \sup_{s \in [0,t]} |\bar{Z}(s)| - M_{\bar{Z}} > u \right\} \leq
		\tilde{\mathbb{Q}}^\xi \left\{ \left| \sup_{s \in [0,t]} |\bar{Z}(s)| - M_{\bar{Z}} \right| > u \right\} \leq  e^{-u^2/2\bar{\sigma}_{\bar{Z}}^2},
	\end{equation}
	where $ M_{\bar{Z}} = \mathbb{E}^{\tilde{\mathbb{Q}}^\xi} \left[ \sup_{s \in [0,t]} |\bar{Z}(s)| \right],$ and $ \bar{\sigma}_{\bar{Z}}^2 = \beta^2 / (2 \bar{\alpha}_0) $.
	Let \( v = u + M_{\bar{Z}} \) with \( v > 2 M_{\bar{Z}} \). Then \eqref{BST1} can be further rearranged as
	\begin{equation}\label{BST2}
		\tilde{\mathbb{Q}}^\xi \left\{ \sup_{s \in [0,t]} |\bar{Z}(s)| > v \right\} \leq  e^{-v^2/8\bar{\sigma}_{\bar{Z}}^2}.
	\end{equation}
	For simplicity, denote \( a_{Y_2} = 2 \overline{k}(\gamma - 1) D_8 t \), and from \eqref{BST2} we get
	\begin{equation*}
		\begin{aligned}
			& \qquad \mathbb{E}^{\tilde{\mathbb{Q}}^\xi} \left[ \exp \left\{ 2 \overline{k}(\gamma - 1) D_8 t \sup_{s \in [0,t]} |\bar{Z}(s)| \right\} \right] = \int_0^\infty \tilde{\mathbb{Q}}^\xi \left( e^{a_{Y_2} \sup_{s \in [0,t]} |\bar{Z}(s)|} > y \right) \textnormal{d}y \\
			& = 1 + \int_0^\infty \tilde{\mathbb{Q}}^\xi \left(  \sup_{s \in [0,t]} |\bar{Z}(s)| > u \right) a_{Y_2} e^{a_{Y_2}u} \textnormal{d}u \\
			& = 1 + \int_0^{2 M_{\bar{Z}}} \tilde{\mathbb{Q}}^\xi \left(  \sup_{s \in [0,t]} |\bar{Z}(s)| > u \right) a_{Y_2} e^{a_{Y_2}u} \textnormal{d}u +
			\int_{2 M_{\bar{Z}}}^{\infty} \tilde{\mathbb{Q}}^\xi \left(  \sup_{s \in [0,t]} |\bar{Z}(s)| > u \right) a_{Y_2} e^{a_{Y_2}u} \textnormal{d}u \\
			& \leq 1 + (e^{2 a_{Y_2} M_{\bar{Z}}} - 1) + 2 \sqrt{2 \pi \bar{\sigma}_{\bar{Z}}^2} a_{Y_2} e^{2 a_{Y_2}^2 \bar{\sigma}_{\bar{Z}}^2}
			= e^{2 a_{Y_2} M_{\bar{Z}}} + 2 \sqrt{2 \pi \bar{\sigma}_{\bar{Z}}^2} a_{Y_2} e^{2 a_{Y_2}^2 \bar{\sigma}_{\bar{Z}}^2} < \infty.
		\end{aligned}
	\end{equation*}
	Therefore,
	\begin{equation*}
		\begin{aligned}
			Y_2(t) & \leq \exp \left\{ 2 \overline{k}(\gamma - 1) D_8 \frac{D_{10}}{\bar{\alpha}_0} t \right\} \mathbb{E}^{\tilde{\mathbb{Q}}^\xi} \left[ \exp \left\{ 2 \overline{k}(\gamma - 1) D_8 \int_0^t |\bar{m}(s)| \, \textnormal{d}s \right\} \right] \\
			& \leq \exp \left\{ 2 \overline{k}(\gamma - 1) D_8 \left(m + \frac{D_{10}}{\bar{\alpha}_0}\right) t \right\} \mathbb{E}^{\tilde{\mathbb{Q}}^\xi} \left[ \exp \left\{ 2 \overline{k}(\gamma - 1) D_8 t \sup_{s \in [0,t]} |\bar{Z}(s)| \right\} \right] < \infty.
		\end{aligned}
	\end{equation*}
	
	Following the same method used to estimate \( Y_2(t) \) above, we can also obtain an estimate for \( Y_1(t) \). We can similarly derive  that
	\begin{equation*}
		m^2(t) \leq \left( \bar{m}(t) + \frac{D_{10}}{\bar{\alpha}_0} \right)^2 \leq 2 \bar{m}^2(t) + \frac{2D_{10}^2}{\bar{\alpha}_0^2},
	\end{equation*}
	\begin{equation*}
		\int_{0}^{t} |\bar{m}^2(s)| \, \textnormal{d}s \leq 2 \int_{0}^{t} |\bar{\mu}^2(s)| \, \textnormal{d}s + \int_{0}^{t} |\bar{Z}^2(s)| \, \textnormal{d}s \leq 2 m^2 t + 2 \left( \sup_{s \in [0,t]} |\bar{Z}(s)| \right)^2 t,
	\end{equation*}
	and
	\begin{equation}\label{BST3}
		\tilde{\mathbb{Q}}^\xi \left\{ \left( \sup_{s \in [0,t]} |\bar{Z}(s)| \right)^2 > v \right\} = \tilde{\mathbb{Q}}^\xi \left\{ \sup_{s \in [0,t]} |\bar{Z}(s)| > \sqrt{v} \right\} \leq  e^{-v/8\bar{\sigma}_{\bar{Z}}^2}.
	\end{equation}
	Then, let \( a_{Y_1} = 8\overline{k}(\gamma - 1) D_7 t \).
	Using \eqref{BST3}, and with condition (H3) ensuring the integrability, we can compute as follows
	\begin{equation*}
		\begin{aligned}
			& \qquad \mathbb{E}^{\tilde{\mathbb{Q}}^\xi} \left[ \exp \left\{ 8\overline{k}(\gamma - 1) D_7 t \left( \sup_{s \in [0,t]} |\bar{Z}(s)| \right)^2 \right\} \right]
			= \int_0^\infty \tilde{\mathbb{Q}}^\xi \left( e^{a_{Y_1} \left( \sup_{s \in [0,t]} |\bar{Z}(s)| \right)^2 } > y \right) \textnormal{d}y \\
			& = 1 + \int_0^\infty \tilde{\mathbb{Q}}^\xi \left( \left( \sup_{s \in [0,t]} |\bar{Z}(s)| \right)^2 > u \right) a_{Y_1} e^{a_{Y_1}u} \textnormal{d}u \\
			& = 1 + \int_0^{2 M_{\bar{Z}}} \tilde{\mathbb{Q}}^\xi \left(  \left( \sup_{s \in [0,t]} |\bar{Z}(s)| \right)^2 > u \right) a_{Y_1} e^{a_{Y_1}u} \textnormal{d}u +
			\int_{2 M_{\bar{Z}}}^{\infty} \tilde{\mathbb{Q}}^\xi \left(  \left( \sup_{s \in [0,t]} |\bar{Z}(s)| \right)^2 > u \right) a_{Y_1} e^{a_{Y_1}u} \textnormal{d}u \\
			& \leq 1 + (e^{2 a_{Y_1} M_{\bar{Z}}} - 1) + \frac{a_{Y_1}}{1/(8 \bar{\sigma}_{\bar{Z}}^2) - a_{Y_1}} e^{- 2 M_{\bar{Z}} \left( 1/(8 \bar{\sigma}_{\bar{Z}}^2) - a_{Y_1} \right) } \\
			& = e^{2 a_{Y_1} M_{\bar{Z}}} + \frac{a_{Y_1}}{1/(8 \bar{\sigma}_{\bar{Z}}^2) - a_{Y_1} } e^{- 2 M_{\bar{Z}} \left( 1/(8 \bar{\sigma}_{\bar{Z}}^2) - a_{Y_1} \right) } < \infty.
		\end{aligned}
	\end{equation*}
	Thus,
	\begin{equation*}
		\begin{aligned}
			Y_1(t) & \leq \exp \left\{ 4\overline{k}(\gamma - 1) D_7 \frac{D_{10}^2}{\bar{\alpha}_0^2} t \right\} \mathbb{E}^{\tilde{\mathbb{Q}}^\xi} \left[ \exp \left\{ 4\overline{k}(\gamma - 1) D_7 \int_0^t |\bar{m}^2(s)| \, \textnormal{d}s \right\} \right] \\
			& \leq \exp \left\{ 4\overline{k}(\gamma - 1) D_7 \left( 2m^2 + \frac{D_{10}^2}{\bar{\alpha}_0^2}\right) t \right\} \mathbb{E}^{\tilde{\mathbb{Q}}^\xi} \left[ \exp \left\{ 8\overline{k}(\gamma - 1) D_7 t \left(\sup_{s \in [0,t]} |\bar{Z}(s)|\right)^2 \right\} \right] < \infty.
		\end{aligned}
	\end{equation*}
	
	Hence we have
	\begin{equation*}
		\sup_{t \in [0, T]} Y_1^{1/2}(t) \cdot Y_2^{1/2}(t) \leq K < \infty,
	\end{equation*}
	where \(K\) denotes a generic constant independent of \(t\), and the same notation \(K\) is used throughout the sequel to represent this meaning.
	Thus, \eqref{EYt} can be further transformed into
	\begin{equation*}
		\mathbb{E}^{\mathbb{Q}^\xi} [X_t^{\overline{k}(1-\gamma)}] \leq K (T - t)^{-\overline{k}(\gamma-1)} \quad \text{for all } t \in [0, T],
	\end{equation*}
	which proves condition (M).

	\medskip
	\noindent \textbf{Step 5: Proof of the Novikov condition used in Step 4.}
	The Novikov condition used in Step 4 is
	\begin{equation*}
		\mathbb{E}^{\mathbb{Q}^\xi} \left[ \exp\left\{ \frac{1}{2} \bar{k}^2(1-\gamma)^2 \int_0^T \left( \sigma^2\tilde{\pi}^2(s) + \lambda\mu_2 \tilde{q}^2(s) \right) \textnormal{d}s \right\} \right] < \infty.
	\end{equation*}
	According to the characterization of the admissible control set \( \mathcal{A}_t \) with respect to \( \tilde{\pi} \) and \( \tilde{q} \), together with \eqref{piqc***}, it suffices to verify
	\begin{equation}\label{Step5-1}
		\mathbb{E}^{\mathbb{Q}^\xi} \left[ \exp\left\{ \frac{\bar{k}^2(1-\gamma)^2}{(\Phi+\gamma)^2} \int_0^T m^2(s) \textnormal{d}s \right\} \right] < \infty.
	\end{equation}
	Note that under the measure \(\mathbb{Q}^\xi\), after substituting \( \xi_1, \xi_2 \) into the measure-transformed process \( m(t) \) \eqref{ROU}, it can be rearranged as
	\begin{equation*}\label{m(t)trans1}
		\textnormal{d} m(t)
		= \left(- \alpha_0 m(t) - \alpha_1 - \alpha_2 \frac{g_m(t,m)}{g(t,m)} \right)\textnormal{d}t + \beta \left[ \rho_1 \textnormal{d} W_{1,t}^{\mathbb{Q}^\xi} + \sqrt{1-\rho_1^2} W_{2,t}^{\mathbb{Q}^\xi} \right],
	\end{equation*}
	where \( \alpha_0 = \alpha + \frac{\Phi\beta\rho_1}{\Phi+\gamma} > 0 \), \( \alpha_1 = \frac{\Phi\beta\rho_1 (a-r)}{(\Phi+\gamma) \sigma} \), \( \alpha_2 = \frac{\Phi k \beta^2 \rho_1^2 (1-\gamma-\Phi)}{(\Phi+\gamma)(1-\gamma)} + \frac{\Phi k \beta^2}{1-\gamma} \).
	
	Similar to the O-U process \eqref{auxiOU} constructed in the analysis of \( Y_2(t) \) in Step 4, the auxiliary O-U process constructed here is given by
	\begin{equation}\label{auxiOU0}
		\textnormal{d}m_0(t) = -\alpha_0 m_0(t) \textnormal{d}t + \beta \left[ \rho_1 \textnormal{d} W_{1,t}^{\mathbb{Q}^\xi} + \sqrt{1-\rho_1^2} \textnormal{d} W_{2,t}^{\mathbb{Q}^\xi} \right], \quad m_0(0) = m.
	\end{equation}
	Define \(\Delta_0(t) = m(t) - m_0(t)\). Then we have
	\begin{equation*}\label{auxiDelta0}
		|\Delta_0(t)| \leq \frac{D_{11}}{\alpha_0}, \ \text{and } m^2(t) \leq 2 m_0^2(t) + \frac{2 D_{11}}{\alpha_0},
	\end{equation*}
	where \( D_{11} = |\alpha_1| + |\alpha_2| K_g \).
	Then, similar to the analysis of \( Y_1(t) \) in Step 4, we can prove that \eqref{Step5-1} holds.

	\medskip
	\noindent \textbf{Step 6: Verification of the condition that \[ \mathbb{E}^{\mathbb{Q}^\xi} \left[\int_t^T |v(s, X_s, m(s)) - J_s| \, \textnormal{d}s \right] < \infty \] for \((\pi, q, c; \xi) \in \mathcal{A}_t\).}
	The proof framework and some technical details of this part are consistent with Step 3 in the proof of Theorem 4.2 of Pu and Zhang \cite{PZ2021}, which is divided into two steps: verifying
	\begin{equation}\label{Step6-1}
		\mathbb{E}^{\mathbb{Q}^\xi} \left[\int_t^T |v(s, X_s, m(s))| \, \textnormal{d}s\right] < \infty,
	\end{equation}
	and
	\begin{equation}\label{Step6-2}
		\mathbb{E}^{\mathbb{Q}^\xi} \left[\int_t^T |J_s| \, \textnormal{d}s\right] < \infty.
	\end{equation}
	The verification of \eqref{Step6-1} follows the same approach as in Pu and Zhang \cite{PZ2021} and is omitted here. For the verification of \eqref{Step6-2}, Pu and Zhang \cite{PZ2021} employed the comparison principle in their derivation. However, since the O-U process involved in this paper does not satisfy the conditions required for the comparison theorem, we handle this part by using the auxiliary O-U process \eqref{auxiOU0} constructed in Step 5. The remaining parts of the verification of \eqref{Step6-2} still follow the same line as in Pu and Zhang \cite{PZ2021}. In summary, the proof of Step 6 is thus completed.
	\qed

	\section*{Disclosure statement}
	
	No potential conflict of interest was reported by the authors.

	\section*{Funding}
	
	This work was supported by the National Natural Science Foundation of China (Grant Nos. 12301603, 12271274).

\end{document}